# Convergence Rates with Inexact Non-expansive Operators

**Jingwei Liang · Jalal Fadili · Gabriel Peyré**

**Abstract** In this paper, we present a convergence rate analysis for the inexact Krasnosel'skiĭ–Mann iteration built from non-expansive operators. The presented results include two main parts: we first establish the global pointwise and ergodic iteration-complexity bounds; then, under a metric sub-regularity assumption, we establish a local linear convergence for the distance of the iterates to the set of fixed points. The obtained results can be applied to analyze the convergence rate of various monotone operator splitting methods in the literature, including the Forward–Backward splitting, the Generalized Forward–Backward, the Douglas–Rachford splitting, alternating direction method of multipliers (ADMM) and Primal–Dual splitting methods. For these methods, we also develop easily verifiable termination criteria for finding an approximate solution, which can be seen as a generalization of the termination criterion for the classical gradient descent method. We finally develop a parallel analysis for the non-stationary Krasnosel'skiĭ–Mann iteration.

## 1 Introduction

### 1.1 Monotone inclusion and operator splitting methods

In various fields of science and engineering, many problems can be cast as solving a structured monotone inclusion problem. A prototypical example that has attracted a wave of interest recently, see e.g. [65,19,58], takes the form

$$\text{Find } x \in \mathcal{H} \text{ such that } 0 \in Bx + \textstyle\sum_{i=1}^n L_i^* \circ A_i \circ L_i x, \tag{1}$$

where $\mathcal{H}$ is a real Hilbert space, $B$ is cocoercive, $A_i$ is a set-valued maximal monotone operator acting on a real Hilbert space $\mathcal{G}_i$, and $L_i : \mathcal{H} \to \mathcal{G}_i$ is a bounded linear operator. Even more complex forms (e.g. with parallel sums (24)) will be discussed in detail in Section 5.

Since the first operator splitting method developed in the 70's for solving structured monotone inclusion problems, the class of splitting methods has been regularly enriched with increasingly sophisticated algorithms as the structure of problems to handle becomes more complex. Splitting methods are iterative algorithms which evaluate (possibly approximately) the individual operators, their resolvents, the linear operators, all separately at various points in the course of iteration, but never the resolvents of sums nor of composition by a linear operator. Popular splitting algorithms to solve special instances of (1) include the Forward–Backward splitting method (FBS) [46,54], the Douglas–Rachford splitting method (DRS) [26,43] and Peaceman–Rachford splitting method (PRS) [55], Alternating Direction method of Multipliers (ADMM) [28,29,30,31]. Other splitting methods were designed to solve (1) or even more complex forms (e.g. with parallel sums), see for instance [15,64,61,12,14,58,19,22,65].

Jingwei Liang, Jalal Fadili
GREYC, CNRS-ENSICAEN-Université de Caen E-mail: {Jingwei.Liang, Jalal.Fadili}@ensicaen.fr

Gabriel Peyré
CNRS and CEREMADE, Université Paris-Dauphine E-mail: gabriel.peyre@ceremade.dauphine.fr

Among the operator splitting methods, many of them (including FBS, DRS, ADMM, and many others) can be cast as the Krasnosel'skiĭ–Mann fixed-point iteration [35,44], possibly in its inexact form to handle errors,

$$z_{k+1} = z_k + \lambda_k(Tz_k + \varepsilon_k - z_k), \tag{2}$$

where $T : \mathcal{H} \to \mathcal{H}$ is a non-expansive operator (Definition 1), and $\mathcal{H}$ is a real Hilbert space $\mathcal{H}$ with norm $\|\cdot\|$. $(\lambda_k)_{k\in\mathbb{N}} \in [0,1]$, and $\varepsilon_k$ is the error of approximating $Tz_k$. The sequence $z_k$ is built in way that it converges (in general weakly [59]) to some fixed point $z^\star$ of $T$, and the latter can be easily related to a solution $x^\star$ of the original problem, e.g. (1).

## 1.2 Contributions

We consider the inexact Krasnosel'skiĭ–Mann iteration (2), and assume that the set of fixed points is non-empty. It is known that a crucial step in proving the convergence of the iterates is to show that $\|z_k - Tz_k\| \to 0$, a property known as asymptotic regularity of $T$ [53]. In this paper, we show that even with errors, the global pointwise and ergodic iteration-complexity bounds, i.e. rates of asymptotic regularity, are respectively $O(1/\sqrt{k})$ and $O(1/k)$. Then under a metric sub-regularity assumption [25], we establish a local linear convergence for the iteration in terms of the distance of the iterates to the set of fixed points. Of course, when the fixed point is unique, the sequence itself converges linearly to it.

Our result can be applied to analyze the convergence behavior of the iterates generated by various monotone operator splitting methods. As stated above, one crucial property of these methods is that they have the equivalent fixed-point formulation (2)[1]. Such methods include the FBS, the Generalized Forward–Backward (GFB) [58], DRS, ADMM, and several Primal–Dual splitting (PDS) methods [65,16]. In particular, for the GFB method developed by two of the co-authors, which addresses the case when $L_i$'s in (1) equal to identity, we demonstrate that $O(1/\epsilon)$ iterations are needed to find a pair $((u_i)_i, g)$ with the termination criterion $\|g + B(\sum_i \omega_i u_i)\|^2 \le \epsilon$, where $g \in \sum_i A_i(u_i)$. This termination criterion can be viewed as a generalization of the classical one based on the norm of the gradient for the gradient descent method [51]. The iteration-complexity improves to $O(1/\sqrt{\epsilon})$ in ergodic sense for the same termination criterion. Similar interpretation is also provided for the DRS/ADMM and PDS considered in Section 5.

We finally study the convergence and rates of the non-stationary version of the Krasnosel'skiĭ–Mann iteration. By absorbing the non-stationarity appropriately into an additional error term, we show that the iteration-complexity bounds developed above remain valid under reasonable conditions. The obtained result is illustrated on the GFB method.

## 1.3 Related work

### 1.3.1 Global iteration-complexity bounds

*Relation with [3,21]* Suppose that $T : \mathcal{C} \to \mathcal{C}$, where $\mathcal{C}$ is a convex *compact* subset of $\mathcal{H}$, and $\varepsilon_k \equiv 0$ (i.e. exact case). Denote $\mathrm{diam}(\mathcal{C})$ the diameter of $\mathcal{C}$. In [3], the authors conjectured that $\|z_k - Tz_k\| \le \mathrm{diam}(\mathcal{C})/\sqrt{\pi \sum_{j=1}^k \lambda_j(1-\lambda_j)}$, and proved it for $\lambda_j \equiv \lambda$. This conjecture was settled in [21, Theorem 1][2]. Our work differs from these in several aspects. For instance, we consider the inexact case without any assumption on boundedness of $\mathcal{C}$. We also establish both the pointwise and ergodic convergence rates, as well as local linear convergence under metric sub-regularity.

However, the extension of the result of [21] to the inexact iteration seems quite intricate. The main reason is that their method of proof relies on the recursive bound in [21, Corollary 3], and exploits some properties of some special functions and an identity for Catalan numbers. These recursions are unfortunately not stable to errors.

*Relation with [33]* The authors in [33] studied the exact version of DRS to solve (1), i.e. $B = 0,\ n = 2$ and $L_i = \mathrm{Id},\ i = 1, 2$. They also supposed the extra assumption that $A_2$ is single-valued. Let

$$e_k = z_k - J_{\gamma A_1}(\mathrm{Id} - \gamma A_2)z_k, \tag{3}$$

[1] In fact, in many cases the fixed-point operator $T$ is even $\alpha$-averaged, see Definition 1.

[2] The authors consider the case where $\mathcal{H}$ is any normed space.

where $J_{\gamma A_1} = (\mathrm{Id} + \gamma A_1)^{-1}$ is the resolvent of $A_1$. Relying on firm non-expansiveness of the resolvent [6], [33] have shown that $\|e_k\| = O(1/\sqrt{k})$.

In fact, it is easy to show that (3) is equivalent to

$$e_k = z_k - \tfrac{1}{2}\big((2J_{\gamma A_1} - \mathrm{Id})(2J_{\gamma A_2} - \mathrm{Id}) + \mathrm{Id}\big)z_k,$$

where single-valuedness of $A_2$ is not needed whatsoever. The operator $\frac{1}{2}((2J_{\gamma A_1} - \mathrm{Id})(2J_{\gamma A_2} - \mathrm{Id}) + \mathrm{Id})$ is firmly non-expansive, and thus fits in our framework. Our results in Section 3 go much beyond this by considering a more general iterative scheme with an operator that is only non-expansive and may be evaluated approximately.

*Relation with HPE* Based on the enlargement of maximal monotone operators, in [62], a hybrid proximal extragradient method (HPE) is introduced to solve monotone inclusion problems of the form

$$\text{Find } x \in \mathcal{H} \ \text{ such that } \ 0 \in Ax.$$

The HPE framework encompasses some splitting algorithms in the literature [48]. The convergence of HPE is established in [62] and in [13] for its inexact version. The pointwise and ergodic iteration-complexities of the exact HPE on a similar error criterion as in our work are established in [48]. Some of the splitting methods we consider in Section 5 are also covered by HPE, hence our iteration-complexity bounds coincide with those of HPE, but only for the ergodic case. While for the pointwise case, our bound is uniform and theirs is not (see further discussion in Remark 3).

### *1.3.2 Local linear convergence*

*Relation with [38, 40]* In [38], local linear convergence of the distance to the set of zeros of a maximal monotone operator using the exact proximal point algorithm (PPA [45, 60]) is established by assuming metric sub-regularity of the operator. Local convergence rate analysis of PPA under a higher-order extension of metric sub-regularity, namely metric $q$-sub-regularity $q \in ]0, 1]$, is conducted in [40]. In our work, metric sub-regularity is assumed on $\mathrm{Id} - T$ with $T$ being the fixed-point operator, i.e. the resolvent, rather than the maximal monotone operator in the case of PPA. Relation between metric sub-regularity of these operators is intricate in general and is beyond the scope of this paper, though we provide an instructive discussion for a simple case at the end of Section 4. Note also that the work of [38, 40] considers PPA only in its classical form, i.e. without errors nor relaxation.

*Relation with [8]* While the first version of this paper was submitted, we became aware of the recent work of [8]. These authors considered (2) with $\lambda_k \equiv \lambda \in [0, 1[$ and $\varepsilon_k \equiv 0$. In [8, Lemma 3.8], they established local linear convergence under bounded linear regularity of $T$ (see [8, Definition 2.1]). The latter coincides exactly with our metric sub-regularity assumption on $\mathrm{Id} - T$ (see (12)). Our result is however more general as it covers the inexact case and $\lambda_k$ is allowed to be iteration-dependent.

*Local linear rate for feasibility problems* In finite dimension, based on strong regularity, [39] proves a local linear convergence of the Method of Alternating Projections (MAP) in the non-convex setting, where the sets are closed and one of which is suitably regular. The linear rate is associated with a modulus of regularity. This is refined later in [7]. In [34], the authors develop a local linear convergence results for the MAP and DRS to solve non-convex feasibility problems, see also [56] for an even more general setting. Their analysis relies on a local version of firm non-expansiveness together with a coercivity condition. It turns out that this coercivity condition holds for mapping $T$ for which the fixed points are isolated and $\mathrm{Id} - T$ is metrically regular. The linear rate they establish, however, imposes a bound on the metric regularity modulus.

*Other local linear rates with DRS* For the case of a sphere intersecting a line, or more generally a proper affine subset, typically in $\mathbb{R}^2$, [9] establishes local linear convergence of DRS. Local linear convergence of (the relaxed) DRS to solve the affine constraint $\ell_1$-minimization problem (basis pursuit) is shown in [23]. For this particular instance, their regularity assumptions can be related to our metric sub-regularity assumption involving the DRS fixed-point operator. Given the level of details this relation requires, we defer it (in an even more general setting) to a forthcoming paper, hence we do not discuss this further.

### 1.4 Paper organization

The organization of the paper is as follows. In Section 2 we introduce our main notations and recall some necessary material on monotone operator theory. Global iteration-complexity bounds and local convergence rate of the inexact Krasnosel'skiĭ–Mann iteration are established in Section 3 and Section 4 respectively. In Section 5, illustrative examples of existing monotone operator splitting methods to which our iteration-complexity result can be applied are described. Extension to the non-stationary case is discussed in Section 6.

## 2 Preliminaries

Throughout the paper, $\mathbb{N}$ is the set of non-negative integers and $\mathcal{H}$ is a real Hilbert space with scalar product $\langle\cdot,\cdot\rangle$, norm $\|\cdot\|$. $\mathrm{Id}$ denotes the identity operator on $\mathcal{H}$. The sub-differential of a proper function $f:\mathcal{H}\to]-\infty,+\infty]$ is the set-valued operator,

$$\partial f:\mathcal{H}\to 2^{\mathcal{H}},\ x\mapsto\big\{g\in\mathcal{H}|(\forall y\in\mathcal{H}),\ \langle y-x,g\rangle+f(x)\leq f(y)\big\}.$$

$\Gamma_0(\mathcal{H})$ denotes the class of proper, lower semi-continuous, convex functions from $\mathcal{H}$ to $]-\infty,+\infty]$. If $f\in\Gamma_0(\mathcal{H})$, then $\mathrm{prox}_f$ denotes the Moreau proximity operator [50], and the Moreau envelope of index $\delta\in]0,+\infty[$ of $f$ is the function,

$${}^{\delta}f:x\mapsto\min_{y\in\mathcal{H}}f(y)+\frac{1}{2\delta}\|y-x\|^2,$$

and its gradient is $\delta^{-1}$-Lipschitz continuous [49].

Let $A:\mathcal{H}\to 2^{\mathcal{H}}$ be a set-valued operator. The domain of $A$ is $\mathrm{dom}A=\{x\in\mathcal{H}|Ax\neq\emptyset\}$, the range of $A$ is $\mathrm{ran}A=\{y\in\mathcal{H}|\exists x\in\mathcal{H}:y\in Ax\}$, the graph of $A$ is the set $\mathrm{gra}A=\{(x,y)\in\mathcal{H}^2|y\in Ax\}$, the inverse of $A$ is the operator whose graph is $\mathrm{gra}A^{-1}=\{(y,x)\in\mathcal{H}^2|x\in A^{-1}y\}$, and its zeros set is $\mathrm{zer}A=\{x\in\mathcal{H}|0\in Ax\}=A^{-1}(0)$.

The resolvent of $A$ is the operator $J_A=(\mathrm{Id}+A)^{-1}$, and the reflection operator associated to $J_A$ is $R_A=2J_A-\mathrm{Id}$.

We denote $\ell^1_+$ the set of summable sequences in $[0,+\infty[$, and define the index set $[\![1,n]\!]=\{1,2,\cdots,n\}$.

A sequence $(x_k)_{k\in\mathbb{N}}$ is said to converge $Q$-linearly to $\tilde{x}$ if there exists a constant $r\in]0,1[$ such that $\frac{\|x_{k+1}-\tilde{x}\|}{\|x_k-\tilde{x}\|}\leq r$, and $(x_k)_{k\in\mathbb{N}}$ is said to converge $R$-linearly to $\tilde{x}$ if $\|x_k-\tilde{x}\|\leq\sigma_k$ and $(\sigma_k)_{k\in\mathbb{N}}$ converges $Q$-linearly to 0.

### 2.1 Non-expansive operators and quasi-Fejér monotone sequences

**Definition 1 (Non-expansive operator)** An operator $T:\mathcal{H}\to\mathcal{H}$ is non-expansive if

$$\forall x,y\in\mathcal{H},\ \|Tx-Ty\|\leq\|x-y\|.$$

For any $\alpha\in]0,1[$, $T$ is $\alpha$-averaged if there exists a non-expansive operator $R$ such that $T=\alpha R+(1-\alpha)\mathrm{Id}$.

We denote $\mathcal{A}(\alpha)$ the class of $\alpha$-averaged operators on $\mathcal{H}$, in particular $\mathcal{A}(\frac{1}{2})$ is the class of firmly non-expansive operators, whose detailed property can be found in [6, Proposition 4.2].

**Lemma 1** *Let $T:\mathcal{H}\to\mathcal{H}$ be an $\alpha$-averaged operator, then* $\frac{1}{2\alpha}(\mathrm{Id}-T)\in\mathcal{A}(\frac{1}{2})$.

*Proof* By definition, there exists a 1-Lipschitz continuous operator $R$ such that $T=\alpha R+(1-\alpha)\mathrm{Id}$, hence

$$\tfrac{1}{2\alpha}(\mathrm{Id}-T)=\tfrac{1}{2\alpha}\big(\mathrm{Id}-(\alpha R+(1-\alpha)\mathrm{Id})\big)=\tfrac{1}{2}\big(\mathrm{Id}+(-R)\big)\in\mathcal{A}(\tfrac{1}{2}).$$

We also collect several useful equivalent characterizations of the firmly non-expansive operators.

**Lemma 2** *The following statements are equivalent:*

(i) *$T$ is firmly non-expansive;*
(ii) $2T-\mathrm{Id}$ *is non-expansive;*
(iii) $\forall x,y\in\mathcal{H},\ \|Tx-Ty\|^2\leq\langle Tx-Ty,x-y\rangle$;
(iv) *$T$ is the resolvent of a maximal monotone operator $A$, i.e. $T=J_A$.*

*Proof* For (i)-(iii), see [6, Proposition 4.2]; For (i) $\Leftrightarrow$ (iv), see [47].

The next lemma shows that the class $\mathcal{A}(\alpha)$ is closed under relaxation, convex combination and composition.

**Lemma 3** *Let $(T_i)_{i\in[\![1,n]\!]}$ be a finite family of non-expansive operators from $\mathcal{H}$ to $\mathcal{H}$, $(\omega_i)_i \in ]0,1]^n$ such that $\sum_i \omega_i = 1$, and let $(\alpha_i)_i \in ]0,1]^n$ such that, for every $i \in [\![1,n]\!]$, $T_i \in \mathcal{A}(\alpha_i)$. Then,*

(i) $(\forall i \in [\![1,n]\!])\ (\forall \lambda_i \in ]0, \frac{1}{\alpha_i}[)$, $\mathrm{Id} + \lambda_i(T_i - \mathrm{Id}) \in \mathcal{A}(\lambda_i \alpha_i)$;
(ii) $T_1 \cdots T_n \in \mathcal{A}(\alpha)$, *with* $\alpha = \frac{n}{n-1+1/(\max_{i\in[\![1,n]\!]}\alpha_i)}$.
(iii) *Let* $\alpha = \max_i\{\alpha_i\}$, *then* $\sum_i \omega_i T_i$ *is* $\alpha$*-averaged.*

*Proof* See [6, Proposition 4.28, 4.30, 4.32].

*Remark 1* For the composite operator $T_1 \cdots T_n$, a sharper bound of $\alpha$ can be obtained for the case $n = 2$ with $\alpha = \frac{\alpha_1+\alpha_2-2\alpha_1\alpha_2}{1-\alpha_1\alpha_2} \in ]0,1[$ ([52, Theorem 3]).

**Definition 2 (Quasi–Fejér monotone sequences)** Let $\Omega$ be a non-empty closed and convex subset of $\mathcal{H}$. A sequence $(x_k)_{k\in\mathbb{N}}$ is *quasi–Fejér monotone* with respect to $\Omega$ if there exists a sequence $(\delta_k)_{k\in\mathbb{N}} \in \ell^1_+$ such that

$$\forall x \in \Omega,\ k \in \mathbb{N}, \quad \|x_{k+1} - x\| \le \|x_k - x\| + \delta_k.$$

## 2.2 Monotone operators

**Definition 3 (Monotone operator)** A set-valued operator $A : \mathcal{H} \to 2^{\mathcal{H}}$ is monotone if

$$\big(\forall (x,u) \in \mathrm{gra}A\big)\big(\forall (y,v) \in \mathrm{gra}A\big),\ \langle x-y, u-v\rangle \ge 0,$$

it is moreover maximally monotone if $\mathrm{gra}A$ is not strictly contained in the graph of any other monotone operator.

**Definition 4 (Cocoercive operator)** An operator $B : \mathcal{H} \to \mathcal{H}$ is called $\beta$-cocoercive for some $\beta \in ]0,+\infty[$ if $\beta B \in \mathcal{A}(\frac{1}{2})$, i.e.,

$$(\forall x, y \in \mathcal{H}),\ \beta\|Bx - By\|^2 \le \langle Bx - By, x - y\rangle.$$

Observe that $\beta$-cocoercivity implies $\frac{1}{\beta}$-Lipschitz continuity.

**Lemma 4** *Let $f : \mathcal{H} \to ]-\infty,+\infty[$ be a convex differentiable function, with $\frac{1}{\beta}$-Lipschitz continuous gradient, $\beta \in ]0,+\infty[$. Then*

(i) $\beta\nabla f \in \mathcal{A}(\frac{1}{2})$, *i.e. is firmly non-expansive;*
(ii) $\mathrm{Id} - \gamma\nabla f \in \mathcal{A}(\frac{\gamma}{2\beta})$ *for* $\gamma \in ]0, 2\beta[$;

*Proof* (i) See [4, Baillon–Haddad theorem]; (ii) See [6, Proposition 4.33].

## 2.3 Product Space

Let $(\omega_i)_i \in ]0,+\infty]^n$, consider $\boldsymbol{\mathcal{H}} = \mathcal{H}^n$ endowed with the scalar product and norm defined as

$$\forall \boldsymbol{x}, \boldsymbol{y} \in \boldsymbol{\mathcal{H}},\ \langle \boldsymbol{x}, \boldsymbol{y}\rangle = \textstyle\sum_{i=1}^n \omega_i \langle x_i, y_i\rangle,\ \|\boldsymbol{x}\| = \sqrt{\textstyle\sum_{i=1}^n \omega_i \|x_i\|^2}.$$

Define $\mathbf{Id}$ the identity operator on $\boldsymbol{\mathcal{H}}$.

Let $\boldsymbol{\mathcal{S}} = \{\boldsymbol{x} = (x_i)_i \in \boldsymbol{\mathcal{H}} | x_1 = \cdots = x_n\}$ and its orthogonal complement $\boldsymbol{\mathcal{S}}^\perp = \big\{\boldsymbol{x} = (x_i)_i \in \boldsymbol{\mathcal{H}} | \sum_{i=1}^n \omega_i x_i = 0\big\} \subset \boldsymbol{\mathcal{H}}$. We also define the canonical isometry $\boldsymbol{C} : \mathcal{H} \to \boldsymbol{\mathcal{S}},\ x \mapsto (x, \cdots, x)$. We have $\forall \boldsymbol{z} \in \boldsymbol{\mathcal{H}}$,

$$P_{\boldsymbol{\mathcal{S}}}(\boldsymbol{z}) = \boldsymbol{C}\big(\textstyle\sum_{i=1}^n \omega_i z_i\big).$$

Clearly $P_{\boldsymbol{\mathcal{S}}}$ is self-adjoint, and its reflection operator is $R_{\boldsymbol{\mathcal{S}}} = 2P_{\boldsymbol{\mathcal{S}}} - \mathbf{Id}$.

Let $\boldsymbol{\gamma} = (\gamma_i)_i \in ]0,+\infty[^n$. For arbitrary maximal monotone operators $A_i, i \in [\![1,n]\!]$ on $\mathcal{H}$, define $\boldsymbol{\gamma A} : \boldsymbol{\mathcal{H}} \to 2^{\boldsymbol{\mathcal{H}}}, \boldsymbol{x} = (x_i)_i \mapsto \times_{i=1}^n \gamma_i A_i x_i$, i.e. its graph is

$$\mathrm{gra}\boldsymbol{\gamma A} = \bigtimes_{i=1}^n \mathrm{gra}\gamma_i A_i = \big\{(\boldsymbol{x}, \boldsymbol{u}) \in \boldsymbol{\mathcal{H}}^2 | \boldsymbol{x} = (x_i)_i,\ \boldsymbol{u} = (u_i)_i,\ u_i \in \gamma_i A_i x_i\big\}.$$

For a single-valued maximal monotone operator $B$, denote $\boldsymbol{B} : \boldsymbol{\mathcal{H}} \to \boldsymbol{\mathcal{H}},\ \boldsymbol{x} = (x_i)_i \mapsto (Bx_i)_i$. It is immediate to check that both $\boldsymbol{\gamma A}$ and $\boldsymbol{B}$ are maximal monotone. We also define the operators $\boldsymbol{B}_{\boldsymbol{\mathcal{S}}} = \boldsymbol{B}P_{\boldsymbol{\mathcal{S}}}$, $J_{\boldsymbol{\gamma A}} = (J_{\gamma_i A_i})_i$ and $R_{\boldsymbol{\gamma A}} = 2J_{\boldsymbol{\gamma A}} - \mathbf{Id}$.

## 3 Global iteration-complexity bounds

In this section, we present the global iteration-complexity bounds of the inexact Krasnosel'skiĭ–Mann iteration.

**Definition 5 (Inexact Krasnosel'skiĭ–Mann iteration)** Let $T : \mathcal{H} \to \mathcal{H}$ be a non-expansive operator such that the set of fixed points $\mathrm{fix}T = \{z \in \mathcal{H} : z = Tz\}$ is non-empty. Let $\lambda_k \in ]0, 1]$, and denote $T_{\lambda_k} = \lambda_k T + (1-\lambda_k)\mathrm{Id}$. Then the inexact Krasnosel'skiĭ–Mann iteration of $T$ is given by

$$z_{k+1} = z_k + \lambda_k(Tz_k + \varepsilon_k - z_k) = T_{\lambda_k} z_k + \lambda_k \varepsilon_k, \tag{4}$$

where $\varepsilon_k$ is the error of approximating $Tz_k$. The error of the iteration is defined as

$$e_k = (\mathrm{Id} - T)z_k = (z_k - z_{k+1})/\lambda_k + \varepsilon_k. \tag{5}$$

Define the following two important notions

$$T' = \mathrm{Id} - T \quad \text{and} \quad \tau_k = \lambda_k(1 - \lambda_k), \tag{6}$$

clearly, $e_k = T'z_k$. We start by collecting some useful properties that characterize the above fixed-point iteration.

**Proposition 1** *The following statements hold,*

(i) $T_{\lambda_k} \in \mathcal{A}(\lambda_k)$ *for* $\lambda_k \in ]0, 1[$*, and if* $T \in \mathcal{A}(\alpha)$*, then* $T_{\lambda_k} \in \mathcal{A}(\lambda_k \alpha)$*;*
(ii) *For any* $z^\star \in \mathrm{fix}T$*,* $z^\star \in \mathrm{fix}T \iff z^\star \in \mathrm{fix}T_{\lambda_k} \iff z^\star \in \mathrm{zer}T'$*;*
(iii) *If* $(\tau_k)_{k\in\mathbb{N}} \notin \ell^1_+$ *and* $(\lambda_k\|\varepsilon_k\|)_{k\in\mathbb{N}} \in \ell^1_+$*, then,*
  (a) $(e_k)_{k\in\mathbb{N}}$ *converges strongly to* $0$*;*
  (b) $(z_k)_{k\in\mathbb{N}}$ *is quasi–Fejér monotone with respect to* $\mathrm{fix}T$*, and converges weakly to a point* $z^\star \in \mathrm{fix}T$*.*

*Proof* (i) Combine Definition 1 and (i) of Lemma 3; (ii) Straightforward; (iii) See [18, Lemma 5.1].

Let's now turn to the properties of $e_k$. Denote $\underline{\tau} = \inf_{k\in\mathbb{N}} \tau_k$, $\overline{\tau} = \sup_{k\in\mathbb{N}} \tau_k$, and $\nu_1 = 2\sup_{k\in\mathbb{N}} \|T_{\lambda_k} z_k - z^\star\| + \sup_{k\in\mathbb{N}} \lambda_k\|\varepsilon_k\|$, $\nu_2 = 2\sup_{k\in\mathbb{N}} \|e_k - e_{k+1}\|$.

**Lemma 5** *For the error term* $e_k$*, the following inequality holds*

$$\tfrac{1}{2\lambda_k}\|e_k - e_{k+1}\|^2 \le \langle e_k - \varepsilon_k, e_k - e_{k+1}\rangle.$$

*Proof* By Lemma 1, $\frac{1}{2}T' \in \mathcal{A}(\frac{1}{2})$. It then follows from Lemma 2(iii) that $\forall p, q \in \mathcal{H}$,

$$\|\tfrac{1}{2}T'(p) - \tfrac{1}{2}T'(q)\|^2 \le \langle p - q, \tfrac{1}{2}T'(p) - \tfrac{1}{2}T'(q)\rangle.$$

Applying this bound with $p = z_k$ and $q = z_{k+1}$, and using the definition of $e_k$ yields the desired result.

**Corollary 1** *If* $T$ *is* $\alpha$*-averaged, the inequality of Lemma 5 becomes*

$$\tfrac{1}{2\alpha\lambda_k}\|e_k - e_{k+1}\|^2 \le \langle e_k - \varepsilon_k, e_k - e_{k+1}\rangle.$$

**Lemma 6** *For* $z^\star \in \mathrm{fix}T$*,* $\lambda_k \in ]0, 1]$*, we have*

$$\|z_{k+1} - z^\star\|^2 \le \|z_k - z^\star\|^2 - \tau_k\|e_k\|^2 + \nu_1\lambda_k\|\varepsilon_k\|.$$

*Proof* By virtue of [6, Corollary 2.14] we get

$$\begin{aligned}
\|z_{k+1} - z^\star\|^2 &= \|T_{\lambda_k} z_k - z^\star + \lambda_k\varepsilon_k\|^2 \\
&\le \|(1-\lambda_k)(z_k - z^\star) + \lambda_k(Tz_k - Tz^\star)\|^2 + \nu_1\lambda_k\|\varepsilon_k\| \\
&= (1-\lambda_k)\|z_k - z^\star\|^2 + \lambda_k\|Tz_k - z^\star\|^2 - \tau_k\|z_k - Tz_k\|^2 + \nu_1\lambda_k\|\varepsilon_k\| \\
&\le \|z_k - z^\star\|^2 - \tau_k\|e_k\|^2 + \nu_1\lambda_k\|\varepsilon_k\|.
\end{aligned}$$

This lemma indicates that $(z_k)_{k\in\mathbb{N}}$ is quasi–Fejér monotone with respect to $\mathrm{fix}T$ as stated in Proposition 1.

**Corollary 2** *If* $T$ *is* $\alpha$*-averaged, the inequality of Lemma 6 holds with* $\lambda_k \in ]0, 1/\alpha]$ *and* $\tau_k = \lambda_k(\frac{1}{\alpha} - \lambda_k)$.

**Lemma 7** *For* $\lambda_k \in ]0,1]$, $(e_k)_{k\in\mathbb{N}}$ *obeys* $\|e_{k+1}\|^2 - \nu_2\|\varepsilon_k\| \le \|e_k\|^2$.

*Proof* We have

$$\begin{aligned}\|e_{k+1}\|^2 &= \|e_{k+1} - e_k + e_k\|^2 = \|e_k\|^2 - 2\langle e_k, e_k - e_{k+1}\rangle + \|e_k - e_{k+1}\|^2 \\ &\le \|e_k\|^2 - 2(1-\lambda_k)\langle e_k - \varepsilon_k, e_k - e_{k+1}\rangle - 2\langle \varepsilon_k, e_k - e_{k+1}\rangle \\ &\le \|e_k\|^2 - \tfrac{1-\lambda_k}{\lambda_k}\|e_k - e_{k+1}\|^2 + \nu_2\|\varepsilon_k\| \le \|e_k\|^2 + \nu_2\|\varepsilon_k\|\end{aligned}$$

where Lemma 5 was used twice in the second and third lines.

We are now in position to the main results of this section. Denote $d_0 = d(z_0, \mathrm{fix}T) = \inf_{z\in\mathrm{fix}T}\|z_0 - z\|$.

**Theorem 1 (Pointwise iteration-complexity bound)** *For the inexact fixed-point iteration (4), if there holds*

$$0 < \inf_{k\in\mathbb{N}} \lambda_k \le \sup_{k\in\mathbb{N}} \lambda_k < 1 \;\text{ and }\; \big((k+1)\|\varepsilon_k\|\big)_{k\in\mathbb{N}} \in \ell_+^1 \,, \tag{7}$$

*then, denoting* $C_1 = \nu_1 \sum_{j\in\mathbb{N}} \lambda_j\|\varepsilon_j\| + \nu_2\overline{\tau}\sum_{\ell\in\mathbb{N}}(\ell+1)\|\varepsilon_\ell\| < +\infty$, *we have*

$$\|e_k\| \le \sqrt{\frac{d_0^2 + C_1}{\underline{\tau}(k+1)}}. \tag{8}$$

*Proof* Condition (7) implies $\underline{\tau} > 0$, $(\tau_k)_{k\in\mathbb{N}} \notin \ell_+^1$ and $(\lambda_k\|\varepsilon_k\|)_{k\in\mathbb{N}} \in \ell_+^1$. Therefore, $z_k$ is quasi-Fejér monotone with respect to $\mathrm{fix}T$ ((iii) of Proposition 1). Thus, $(\|e_k\|)_{k\in\mathbb{N}}$ and $(\|z_k - z^\star\|)_{k\in\mathbb{N}}$ are bounded for any $z^\star \in \mathrm{fix}T$. Hence $\nu_1, \nu_2$ and $C_1$ are bounded constants. Choose $z^\star$ such that $d_0 = \|z_0 - z^\star\|$. From Lemma 6, $\forall k \in \mathbb{N}$,

$$\tau_k\|e_k\|^2 \le \|z_k - z^\star\|^2 - \|z_{k+1} - z^\star\|^2 + \nu_1\lambda_k\|\varepsilon_k\|.$$

Summing up from $j = 0$ to $k$,

$$\textstyle\sum_{j=0}^k \tau_j\|e_j\|^2 \le \|z_0 - z^\star\|^2 - \|z_{k+1} - z^\star\|^2 + \nu_1\sum_{j=0}^k \lambda_j\|\varepsilon_j\|. \tag{9}$$

From Lemma 7, we have $\forall j \le k$,

$$\textstyle\|e_k\|^2 - \nu_2\sum_{\ell=j}^{k-1}\|\varepsilon_\ell\| \le \|e_j\|^2.$$

Substituting this back into (9) yields,

$$\begin{aligned}\textstyle\big(\sum_{j=0}^k\tau_j\big)\|e_k\|^2 &\le \textstyle\sum_{j=0}^k\tau_j\|e_j\|^2 + \nu_2\sum_{j=0}^k\tau_j\sum_{\ell=j}^{k-1}\|\varepsilon_\ell\| \\ &\le \textstyle d_0^2 + \nu_1\sum_{j=0}^k\lambda_j\|\varepsilon_j\| + \nu_2\sum_{j=0}^k\tau_j\sum_{\ell=j}^{k-1}\|\varepsilon_\ell\|.\end{aligned} \tag{10}$$

Finally, since $(k+1)\underline{\tau} \le \sum_{j=0}^k\tau_j$, we get,

$$\textstyle(k+1)\underline{\tau}\|e_k\|^2 \le d_0^2 + \nu_1\sum_{j=0}^k\lambda_j\|\varepsilon_j\| + \nu_2\overline{\tau}\sum_{\ell=0}^{k-1}(\ell+1)\|\varepsilon_\ell\|,$$

which leads to the desired result (8).

*Remark 2*

(i) If $T$ is $\alpha$-averaged, then condition (7) of Theorem 1 on $\lambda_k$ changes to $0 < \inf_{k\in\mathbb{N}}\lambda_k \le \sup_{k\in\mathbb{N}}\lambda_k < \frac{1}{\alpha}$.
(ii) Since finding $z^\star \in \mathrm{fix}T$ is equivalent to finding a zero of $T'$ ((ii) of Proposition 1), Theorem 1 tells us that $O(1/\epsilon)$ iterations are needed for (4) to reach an $\epsilon$-accurate in terms of the error criterion $\|T'z_k\|^2 \le \epsilon$.
(iii) For the case of first-order methods for solving smooth optimization problems, i.e. the gradient descent where $T'$ is just the gradient, the obtained pointwise bound is the best-known complexity bound [51].
(iv) When the fixed-point iteration (4) is exact, the sequence $(\|e_k\|)_{k\in\mathbb{N}}$ is non-increasing (Lemma 7), hence we get $\|e_k\| \le d_0/\sqrt{\sum_{j=0}^k\tau_j}$, which recovers the result of [21, Proposition 11]. Note that we provide a sharper monotonicity property compared to them.

We now turn to the ergodic iteration-complexity bound of (4). For this, let's define $\Lambda_k = \sum_{j=0}^k \lambda_j$ and $\bar{e}_k = \frac{1}{\Lambda_k}\sum_{j=0}^k \lambda_j e_j$.

**Theorem 2 (Ergodic iteration-complexity bound)** *Suppose that* $C_2 = \sum_{k\in\mathbb{N}} \lambda_k\|\varepsilon_k\| < +\infty$. *Then,*

$$\|\bar{e}_k\| \le \frac{2(d_0 + C_2)}{\Lambda_k}.$$

*In particular, if* $\inf_{k\in\mathbb{N}} \lambda_k > 0$, *then* $\|\bar{e}_k\| = O(1/k)$.

*Proof* Again, let $z^\star \in \mathrm{fix}T$ such that $d_0 = \|z_0 - z^\star\|$. Since $T_{\lambda_k}$ is non-expansive, we have

$$\begin{aligned}\|z_{k+1} - z^\star\| &= \|T_{\lambda_k} z_k - T_{\lambda_k} z^\star + \lambda_k\varepsilon_k\| \le \|z_k - z^\star\| + \lambda_k\|\varepsilon_k\| \\ &\le \|z_{k-1} - z^\star\| + \textstyle\sum_{j=k-1}^{k}\lambda_j\|\varepsilon_j\| \le \|z_0 - z^\star\| + \textstyle\sum_{j=0}^{k}\lambda_j\|\varepsilon_j\|.\end{aligned}$$

This together with the definition of $\bar{e}_k$ yields

$$\begin{aligned}\|\bar{e}_k\| &= \|\tfrac{1}{\Lambda_k}\textstyle\sum_{j=0}^{k} \lambda_j e_j\| = \tfrac{1}{\Lambda_k}\|\textstyle\sum_{j=0}^{k} (z_j - z_{j+1}) + \textstyle\sum_{j=0}^{k}\lambda_j\varepsilon_j\| \\ &\le \tfrac{1}{\Lambda_k}\big(\|z_0 - z^\star\| + \|z_{k+1} - z^\star\| + \textstyle\sum_{j=0}^{k}\lambda_j\|\varepsilon_j\|\big) \le \tfrac{2(d_0+C_2)}{\Lambda_k}.\end{aligned}$$

Again, this result holds when $T$ is $\alpha$-averaged, where now $\lambda_k$ is allowed to vary in $]0, 1/\alpha]$.

*Remark 3*

(i) As in the pointwise case, Theorem 2 holds when $T$ is $\alpha$-averaged, where now $\lambda_k$ is allowed to vary in $]0, 1/\alpha]$.

(ii) When $T$ is firmly non-expansive, i.e. the resolvent of a maximal monotone operator ((iv) of Lemma 2), an $O(1/k)$ ergodic convergence rate is also established in [48] with summable enlargement errors. For the methods which can also be cast in the HPE framework, our result coincides with the one in [48]. Note that in that work, they handled the non-stationary case (i.e. the parameter of the resolvent varies); see also our extension to the non-stationary case in Section 6. For the case without errors, we recover also the bound in [2].

From Theorem 1 and 2, it is immediate to get the convergence rate bounds on the sequence $(\|z_k - z_{k+1}\|)_{k\in\mathbb{N}}$ in the exact case. To lighten the notation, let $v_k = z_k - z_{k+1}$ and $\bar{v}_k = \frac{1}{k+1}\sum_{j=0}^{k} v_j$.

**Corollary 3** *Assume that* $\varepsilon_k = 0$ *for all* $k \in \mathbb{N}$.

(i) *If* $0 < \inf_{k\in\mathbb{N}} \lambda_k \le \sup_{k\in\mathbb{N}} \lambda_k < 1$, *then* $\|v_k\| \le \frac{d_0}{\sqrt{\underline{\tau}(k+1)}}$;

(ii) *If* $\underline{\lambda} = \inf_{k\in\mathbb{N}} \lambda_k > 0$, *then* $\|\bar{v}^k\| \le \frac{2d_0}{k+1}$.

*Proof*

(i) By definition $v_k = \lambda_k e_k$, then from (8) we have

$$\|v^k\| = \|\lambda_k e_k\| \le \lambda_k \frac{d_0}{\sqrt{\underline{\tau}(k+1)}} \le \frac{d_0}{\sqrt{\underline{\tau}(k+1)}}.$$

(ii) A direct result of Theorem 2 by replacing $\Lambda_k$ with $k + 1$.

## 4 Local linear convergence rate

In the literature, for many splitting algorithms applied to a range of optimization problems, the following typical convergence profile has been observed in practice. Globally the algorithm converges sub-linearly, and after a sufficiently large number of iterations, the algorithm enters a new regime where a local linear convergence takes over. This has been for instance observed (and sometimes proved) for DRS or FBS when solving sparsity-enforcing minimization problems, see e.g. [23,41].

In this section, we study the rationale underlying this local linear convergence behavior. Our analysis relies on a metric sub-regularity assumption of $T' = \mathrm{Id} - T$ (Definition 5).

**Definition 6 (Metric sub-regularity [25])** A set-valued mapping $F : \mathcal{H} \to 2^{\mathcal{H}}$ is called metrically sub-regular at $\tilde{z}$ for $\tilde{u} \in F(\tilde{z})$ if there exists $\kappa \ge 0$ along with neighborhood $\mathcal{Z}$ of $\tilde{z}$ such that

$$d(z, F^{-1}\tilde{u}) \le \kappa\, d(\tilde{u}, Fz), \quad \forall z \in \mathcal{Z}. \tag{11}$$

The infimum of $\kappa$ for which (11) holds is the modulus of metric sub-regularity, denoted by $\mathrm{subreg}(F; \tilde{z}|\tilde{u})$. The absence of metric regularity is signaled by $\mathrm{subreg}(F; \tilde{z}|\tilde{u}) = +\infty$.

Metric sub-regularity implies that, for any $z \in \mathcal{Z}$, $d(\tilde{u}, Fz)$ is bounded below. The metric (sub)regularity of multifunctions plays a crucial role in modern variational analysis and optimization. These properties are a key to study the stability of solutions of generalized equations, see the dedicated monograph [25].

Let's specialize this notion to $T'$ and $\tilde{u} = 0$. Since $T'$ is single-valued and $\mathrm{zer}T' = \mathrm{fix}T$, from (11), metric sub-regularity of $T'$ at some $z^\star \in \mathrm{fix}T$ for 0 is equivalent to

$$d(z, \mathrm{fix}T) \leq \kappa \|T'z\|, \quad \forall z \in \mathcal{Z}. \tag{12}$$

There are several concrete examples of operators $T$ where $T'$ fulfills (12).

*Example 1 (Projector)* Suppose that $T = P_{\mathcal{C}}$, for $\mathcal{C}$ a non-empty closed convex subset of $\mathcal{H}$. Then $\mathrm{fix}T = \mathcal{C}$ and $d(z, \mathcal{C}) = \|z - P_{\mathcal{C}} z\| = \|T'z\|$. Thus $T'$ is metrically sub-regular at any $z^\star \in \mathcal{C}$ for 0 with $\mathcal{Z} = \mathcal{H}$ and modulus 1.

Using the relation between metric sub-regularity of $T'$ and bounded linear regularity of $T$ as defined in [8] (see our discussion in Section 1.3.2), other examples can be deduced for instance from [8, Example 2.3 and 2.5]. Two other instructive examples, one on DRS with two subspaces and the second on gradient descent will be discussed at the end of the section.

Another interesting situation is when $T$ is firmly non-expansive, so that $T = J_F$ for some maximal monotone operator $F$ ((iv) of Lemma 2), in which case (4) is the relaxed inexact PPA. Let $z^\star \in \mathrm{fix}T = \mathrm{zer}F$ and suppose that $0 \in \mathrm{zer}F$. If $F$ is metrically sub-regular at $z^\star$ for 0 with modulus $\gamma$, then

$$d(z, \mathrm{zer}F) = d(z, \mathrm{fix}T) \leq \gamma d(0, Fz) = \gamma d(0, T^{-1}z - z) = \gamma \inf_{v \in T^{-1}z} \|v - z\|, \ \forall z \in \mathcal{Z}.$$

Thus for all $w$ such that $z = Tw \in \mathcal{Z}$, applying the previous inequality and using the triangle inequality, we get

$$d(Tw, \mathrm{fix}T) \leq \gamma \|T'w\| \quad \text{and} \quad d(w, \mathrm{fix}T) \leq (1 + \gamma)\|T'w\|.$$

Clearly, this is closely related, though not equivalent, to metric sub-regularity of $T'$.

Metric sub-regularity implies that (12) gives an estimate for how far a point $z$ is from being the fixed-point set of $T$ in terms of the residual $\|z - Tz\|$. This is the rationale behind using such a regularity assumption on the operator $T'$ to quantify the convergence rate on $d(z_k, \mathrm{fix}T)$. Thus, starting from $z_0 \in \mathcal{H}$, and by virtue of Theorem 1, one can recover a $O(1/\sqrt{k})$ rate on $d(z_k, \mathrm{fix}T)$. In fact, we can do even better as shown in the following theorem. We use the shorthand notation $d_k = d(z_k, \mathrm{fix}T)$.

**Theorem 3 (Local convergence rate)** *Let $z^\star \in \mathrm{fix}T$, assume $T'$ is metrically sub-regular at $z^\star$ with neighborhood $\mathcal{Z}$ of $z^\star$, let $\kappa > \mathrm{subreg}(T'; z^\star|0)$, $\lambda_k \in ]0, 1]$. Given a ball $\mathbb{B}_a(z^\star) \subset \mathcal{Z}, a \geq 0$, suppose $C_2 = \sum_{k \in \mathbb{N}} \lambda_k \|\varepsilon_k\|$ is small enough such that*

$$\mathbb{B}_{(a+C_2)}(z^\star) \subseteq \mathcal{Z}.$$

*Then for any starting point $z_0 \in \mathbb{B}_a(z^\star)$, we have for all $k \in \mathbb{N}$,*

$$d_{k+1}^2 \leq \zeta_k \, d_k^2 + c_k, \quad \text{where} \quad \zeta_k = \begin{cases} 1 - \frac{\tau_k}{\kappa^2}, & \text{if } \tau_k/\kappa^2 \in ]0, 1] \\ \frac{\kappa^2}{\kappa^2 + \tau_k}, & \text{otherwise} \end{cases} \in [0, 1[, \tag{13}$$

*and $c_k = \nu_1 \lambda_k \|\varepsilon_k\|$. Moreover,*

(i) *$d_k$ converges to 0 if $(\tau_k)_{k \in \mathbb{N}} \notin \ell_+^1$.*

(ii) *Let $\chi_k = \prod_{j=0}^k \zeta_j$, if $\chi = \limsup_{k \to +\infty} \sqrt[k]{\chi_k} < 1$, then $(d_k^2)_{k \in \mathbb{N}} \in \ell_+^1$. When $\varepsilon_k = 0$, then $\lim_{k \to +\infty} \sqrt[k]{d_k} < 1$, which is R-linear convergence.*

(iii) *If $0 < \inf_{k \in \mathbb{N}} \lambda_k \leq \sup_{k \in \mathbb{N}} \lambda_k < 1$, then there exists $\zeta \in (0, 1)$ such that*

$$d_{k+1}^2 \leq \zeta^k \big(d_0^2 + \textstyle\sum_{j=0}^k \zeta^{-j+1} c_j\big).$$

*Proof* Let constants $b > a > 0$ such that $T'$ is metrically sub-regular at $z^\star$ with modulus $\kappa$. Make the radius $a$ smaller if necessary so that $\mathbb{B}_a(z^\star) \subset \mathbb{B}_b(z^\star) \subseteq \mathcal{Z}$ and

$$a + C_2 \leq b.$$

Pick $z_0 \in \mathbb{B}_a(z^\star)$. If $z_0 = z^\star$ then, take $z_k = z^\star$ for all $k \in \mathbb{N}$ and there is nothing more to prove. If not, then from Lemma 6 that for any $\tilde{z} \in \mathrm{fix}T$,

$$\|z_{k+1} - \tilde{z}\| \le \|z_k - \tilde{z}\| + \lambda_k\|\varepsilon_k\| \le \cdots \le \|z_0 - \tilde{z}\| + \textstyle\sum_{j=0}^{k}\lambda_j\|\varepsilon_j\| \le a + C_2 \le b,$$

which implies that starting from any point $z_0 \in \mathbb{B}_a(z^\star)$, $z_k \in \mathbb{B}_b(z^\star)$ holds for all $k \in \mathbb{N}$. Now for $\forall k \in \mathbb{N}$, let $\tilde{z} \in \mathrm{fix}T$ be such that $d_k = \|z_k - \tilde{z}\|$ and $c_k = \nu_1\lambda_k\|\varepsilon_k\|$, then by virtue of the metric sub-regularity of $T'$ and Lemma 6, we have

$$d_{k+1}^2 \le \|z_{k+1} - \tilde{z}\|^2 \le \|z_k - \tilde{z}\|^2 - \tau_k\|T'z_k - T'\tilde{z}\|^2 + c_k \le d_k^2 - \tfrac{\tau_k}{\kappa^2}d_k^2 + c_k \tag{14}$$

$$\le d_k^2 - \tfrac{\tau_k}{\kappa^2}(d_{k+1}^2 - c_k) + c_k = d_k^2 - \tfrac{\tau_k}{\kappa^2}d_{k+1}^2 + \tfrac{1}{\zeta_k}c_k. \tag{15}$$

If $\tau_k/\kappa^2 \in ]0,1[$, then from (14) we have $d_{k+1}^2 \le (1 - \frac{\tau_k}{\kappa^2})d_k^2 + c_k$, or if $1 \le \tau_k/\kappa^2$, (15) produces $d_{k+1}^2 \le \frac{\kappa^2}{\kappa^2+\tau_k}d_k^2 + c_k$, $\lambda_k \in ]0,1]$ ensures $\kappa^2/(\kappa^2+\tau_k) \in ]0,1]$. Therefore, we have

$$\zeta_k = \begin{cases} 1 - \frac{\tau_k}{\kappa^2}, & \text{if } \tau_k/\kappa^2 \in ]0,1[ \\ \frac{\kappa^2}{\kappa^2+\tau_k}, & \text{if } 1 \le \tau_k/\kappa^2 \end{cases} \in ]0,1].$$

Furthermore,

$$d_{k+1}^2 \le \zeta_k\, d_k^2 + c_k \le \cdots \le \chi_k\, d_0^2 + \textstyle\sum_{j=0}^{k}\phi_{k-j}c_j \le \chi_k\, d_0^2 + \sum_{j=0}^{k}c_j, \tag{16}$$

where $\chi_k = \prod_{j=0}^{k}\zeta_j$ and $\phi_{k-j} = \prod_{\ell=j+1}^{k}\zeta_\ell$.

(i) From (16) we have $d_{k+1}^2 \le d_k^2 + c_k$, then the $d_k^2 \to d \ge 0$ ([57, Lemma 2.2.2]). If $(\tau_k)_{k\in\mathbb{N}} \notin \ell_+^1$, then $\|e_k\| \to 0$ (Theorem 1), and by metric sub-regularity we have $d_k \le \kappa\|e_k\|$, therefore $d = 0$;

(ii) If $\chi = \limsup_{k\to+\infty}\sqrt[k]{\chi_k} < 1$, then $\lim_{k\to+\infty}\chi_k = 0$ and $(\chi_k)_{k\in\mathbb{N}}$, $(\phi_k)_{k\in\mathbb{N}} \in \ell_+^1$. Since also $(c_k)_{k\in\mathbb{N}} \in \ell_+^1$, hence the convolution $\big(\sum_{j=0}^{k}\phi_{k-j}c_j\big)_{k\in\mathbb{N}} \in \ell_+^1$, as a result, $\big(\chi_k\, d_0^2 + \sum_{j=0}^{k}\phi_{k-j}c_j\big)_{k\in\mathbb{N}} \in \ell_+^1$ and so is $(d_k^2)_{k\in\mathbb{N}}$. If $\varepsilon_k = 0$, then from (16) we have $\lim_{k\to+\infty}\sqrt[k]{d_k} \le \limsup_{k\to+\infty}\sqrt[k]{\chi_k} < 1$.

(iii) If $0 < \inf_{k\in\mathbb{N}}\lambda_k \le \sup_{k\in\mathbb{N}}\lambda_k < 1$, then there exists $\zeta = \sup_{k\in\mathbb{N}}\kappa^2/(\kappa^2+\tau_k) < 1$ which concludes the result.

*Remark 4*

(i) When the fixed point is a singleton, Theorem 3 holds replacing $d_k$ by $\|z_k - z^\star\|$.

(ii) For simplicity, suppose the iteration is exact, and let $z^\star \in \mathrm{fix}T$ such that $d_k = \|z_k - z^\star\|$. Then we have

$$\|e_k\|^2 = \|z_k - z^\star + Tz^\star - Tz_k\|^2 \le 4d_k^2,$$

which means locally, $\|e_k\|$ also converges linearly to 0 given that $(\tau_k)_{k\in\mathbb{N}} \notin \ell_+^1$.

(iii) As far as the claim in (iii) of Theorem 3 is concerned, if $\exists\xi \in ]0,1]$ such that $c_k = O(\xi^k)$, then
   (a) If $\xi < \zeta$, then $d_{k+1}^2 = O(\zeta^k)$;
   (b) If $\xi = \zeta$, then $d_{k+1}^2 = O(k\zeta^k) = o\big((\zeta+\delta)^k\big)$, $\delta > 0$;
   (c) If $\xi > \zeta$, then $d_{k+1}^2 = o(\xi^k)$.

Theorem 3 extends readily to the $\alpha$-averaged case.

**Corollary 4** *If $T$ is $\alpha$-averaged, then Theorem 3 holds substituting $\lambda_k\alpha$ for $\lambda_k$ and $\kappa\alpha$ for $\kappa$.*

*Remark 5*

(i) When $\lambda_k \equiv \lambda$ and $\varepsilon_k \equiv 0$, our second rate estimate in (13) encompasses that of [8, Lemma 3.8].

(ii) Equivalent characterizations of metric sub-regularity can be given, for instance in terms of derivative criteria. In particular, as $T'$ is single-valued, metric sub-regularity of $T'$ holds if $T'$ is differentiable on a neighborhood of $z^\star$ with non-singular derivatives at $z$ around $z^\star$, and the operator norms of their inverses are uniformly bounded [24, Theorem 1.2]. Computing the metric regularity modulus $\kappa$ is however far from obvious in general even for the differentiable case.

**Examples**

*DRS with two subspaces in* $\mathbb{R}^2$ Let $\mathcal{U}$ and $\mathcal{V}$ two subspaces in $\mathbb{R}^2$ forming an angle of $\theta \in ]0, \pi/2]$. Consider the problem of finding $\mathcal{U} \cap \mathcal{V} = \{0\}$ using DRS. It is well-known the DRS fixed point operator $T_{\text{DRS}}$ is firmly non-expansive. Moreover, following [8, Example 2.3 and 2.5], it can be shown $T'_{\text{DRS}} = \text{Id} - T_{\text{DRS}}$ is metrically sub-regular with modulus $1/\sin\theta$. It then follows from (13) that the rate estimate is $\zeta_k = 1 - (2-\lambda_k)\lambda_k \sin^2\theta = (1-\lambda_k)^2 + (2-\lambda_k)\lambda_k \cos^2\theta \in ]0,1[$, for $\lambda_k \in ]0,2[$. This is exactly the optimal rate estimate provided in [42,5]. Observe that as remarked in [42], the best rate $\zeta_k \equiv \cos^2\theta$ is obtained for $\lambda_k \equiv 1$ (no relaxation).

*Gradient descent* We now consider the example of minimizing a convex function $f : \mathbb{R}^n \to \mathbb{R}$ which is locally $C^2$ and strongly convex on $\Omega \subset \mathbb{R}^n$ and has Lipschitz gradient. Clearly $f$ has a unique minimizer, and without losing generality, we denote it as 0. Define $\delta_{\text{m}}, \delta_{\text{M}}$ the local strong convexity and Lipschitz modulus of $f$ respectively, then $\frac{1}{\delta_{\text{m}}}$ is equivalent to metric sub-regularity of $\nabla f$ at 0 for 0 [1, Theorem 3.5].

For simplicity, consider the non-relaxed gradient descent method [57] for minimizing $f$ with constant step-size, i.e.

$$x_{k+1} = x_k - \gamma \nabla f(x_k),$$

where $\gamma \in ]0, \frac{2}{\delta_{\text{M}}}[$. This can be cast in our above framework by setting $T = \text{Id} - \gamma\nabla f$ which is $\frac{\gamma\delta_{\text{M}}}{2}$-averaged (Lemma 4). Hence $T' \stackrel{\text{def}}{=} \text{Id} - T = \gamma\nabla f$ which is continuously differentiable on any neighborhood $\mathcal{Z} \subset \Omega$. For any $x \in \mathcal{Z}$, the Jacobian of $T'$ is just $\gamma\nabla^2 f(x)$, where $\nabla^2 f(x)$ is the Hessian of $f$ at $x$, which is non-singular and its inverse is uniformly bounded by $\frac{1}{\delta_{\text{m}}}$ from local strong convexity. Thus, by virtue of [25, Theorem 4B.1], $T'$ is metrically regular, hence sub-regular, and the metric regularity modulus $\kappa$ is precisely $\frac{1}{\gamma\delta_{\text{m}}}$. In fact we could have anticipated this directly from the local strong monotonicity of $\nabla f$. Specializing the rate of Theorem 3, we get

$$\zeta = 1 - \frac{t(2-t)}{\text{cnd}^2} \in [0,1[,$$

where we set $t = \gamma\delta_{\text{M}} \in ]0,2[$, and $\text{cnd} = \delta_{\text{M}}/\delta_{\text{m}}$ can be seen as the condition number of the Hessian of $f$ along $\Omega$. It is obvious that the smallest $\zeta$[3] is attained for $t = 1$, i.e. $\gamma = 1/\delta_{\text{M}}$.

The observed and theoretical convergence profiles of $\|e_k\| = \|x_k - x_{k+1}\|$ are illustrated in Figure 1 where gradient descent was run with $\gamma \in \{\frac{1}{2\delta_{\text{M}}}, \frac{1}{\delta_{\text{M}}}\}$. As predicted by our result, the convergence profile exhibits two regimes, a global sub-linear one, and then a local linear one.

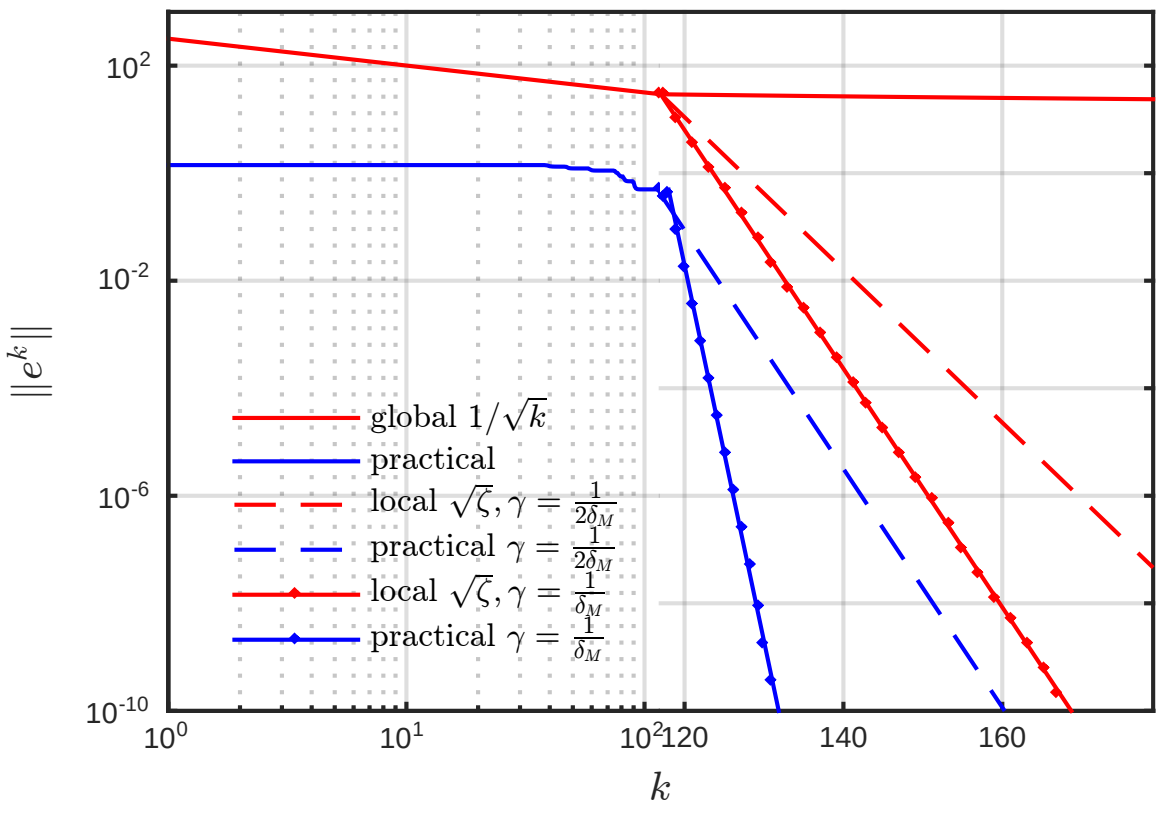

**Fig. 1** Global and local convergence profiles of gradient descent to minimize $f$ with $(\delta_{\text{m}}, \delta_{\text{M}}) = (0.8, 1)$. For $\gamma = \frac{1}{2\delta_{\text{M}}}$, the observed and theoretical rates are respectively 0.60 and 0.72. For $\gamma = \frac{1}{\delta_{\text{M}}}$, they are 0.20 and 0.60 respectively.

[3] One may observe that $\gamma$ can be modified locally, once the iterates enter the appropriate neighborhood, to get the usual optimal linear rate of gradient descent with a strongly convex objective [57]. But this is not our aim here.

## 5 Applications

In this section, we apply the obtained results to conduct quantitative convergence analysis of a class of monotone operator splitting methods in the literature, and mainly focus on the global iteration-complexity bounds. As stated in the introduction, we will rely on the fact that all the considered iterative schemes (GFB/FBS, DRS/ADMM and PDS) can be cast as Krasnosel'skiĭ–Mann iteration. Furthermore, based on the structure of the corresponding monotone inclusion problem, we also derive specific error criteria which can serve as termination tests.

### 5.1 Generalized Forward–Backward splitting

The GFB algorithm [58] addresses the monotone inclusion (1) when $L_i = \mathrm{Id}, i \in [\![1,n]\!]$, and $B$ is $\beta$-cocoercive. By lifting the problem into product space, GFB achieves full splitting by applying the $A_i$'s implicitly and $B$ explicitly. Let $(\omega_i)_i \in ]0,1[^n$ such that $\sum_{i=1}^n \omega_i = 1$, $\gamma \in ]0, 2\beta[$, $\lambda_k \in ]0, \frac{4\beta-\gamma}{2\beta}[$, $(b_k)_{k\in\mathbb{N}}$ and $(a_{i,k})_{k\in\mathbb{N}}$ are error terms. The GFB iteration reads

$$\left\lfloor \begin{array}{l} \text{For } i = 1, \cdots, n \\ \left\lfloor z_{i,k+1} = z_{i,k} + \lambda_k\big(J_{\frac{\gamma}{\omega_i}A_i}(2x_k - z_{i,k} - \gamma B x_k + b_k) + a_{i,k} - x_k\big), \right. \\ x_{k+1} = \sum_{i=1}^n \omega_i z_{i,k+1}. \end{array} \right. \tag{17}$$

GFB recovers FBS for $n = 1$, and when $B = 0$, GFB recovers DRS in the *product space*. In the literature, the convergence properties of the FBS and DRS have been extensively studied, see [6] and references therein. In this section, we mainly focus on the GFB method.

Recall the notation introduced in Section 2.3. Denote $u_{i,k+1} = J_{\frac{\gamma}{\omega_i}A_i}(2x_k - z_{i,k} - \gamma B x_k),\ i \in [\![1,n]\!]$, then (17) equivalently reads

$$z_{i,k+1} = z_{i,k} + \lambda_k(u_{i,k+1} + \varepsilon_{i,k} - x_k),$$

where $\varepsilon_{i,k} = \big(J_{\frac{\gamma}{\omega_i}A_i}(2x_k - z_{i,k} - \gamma B x_k + b_k) + a_{i,k}\big) - u_{i,k+1}$. Let $\boldsymbol{\varepsilon}_{1,k} = \boldsymbol{C}(b_k)$, $\boldsymbol{\varepsilon}_{2,k} = (a_{i,k})_i$, and $\boldsymbol{T}_{1,\gamma} = \frac{1}{2}(R_{\gamma\boldsymbol{A}} R_{\boldsymbol{\mathcal{S}}} + \mathbf{Id})$, $\boldsymbol{T}_{2,\gamma} = \mathbf{Id} - \gamma \boldsymbol{B}_{\boldsymbol{\mathcal{S}}}$.

Using Remark 1, the composed operator $\boldsymbol{T}_\gamma = \boldsymbol{T}_{1,\gamma} \circ \boldsymbol{T}_{2,\gamma}$ is $\frac{2\beta}{4\beta-\gamma}$-averaged. Moreover, it was shown in [58] that the fixed-point iteration of GFB reads,

$$\boldsymbol{z}_{k+1} = \boldsymbol{z}_k + \lambda_k\big(\boldsymbol{T}_{1,\gamma}(\boldsymbol{T}_{2,\gamma}\boldsymbol{z}_k + \boldsymbol{\varepsilon}_{1,k}) + \boldsymbol{\varepsilon}_{2,k} - \boldsymbol{z}_k\big) = \boldsymbol{z}_k + \lambda_k(\boldsymbol{T}_\gamma \boldsymbol{z}_k + \boldsymbol{\varepsilon}_k - \boldsymbol{z}_k), \tag{18}$$

where $\boldsymbol{\varepsilon}_k = \big(\boldsymbol{T}_{1,\gamma}(\boldsymbol{T}_{2,\gamma}\boldsymbol{z}_k + \boldsymbol{\varepsilon}_{1,k}) + \boldsymbol{\varepsilon}_{2,k}\big) - \boldsymbol{T}_\gamma \boldsymbol{z}_k$.

Obviously, (18) is exactly in the form of (4). Therefore, the GFB iterates converge weakly [58, Theorem 4.1], and obey the iteration-complexity bounds in Theorem 1 and 2. Moreover, we can establish certain convergence rates for the structured monotone inclusion (1).

Let $g_{k+1} = \frac{1}{\gamma}x_k - Bx_k - \frac{1}{\gamma}\sum_i \omega_i u_{i,k+1}$. Define $\bar{u}_{i,k+1} = \frac{1}{k+1}\sum_{j=0}^k u_{i,j+1},\ i \in [\![1,n]\!]$, $\bar{x}_k = \frac{1}{k+1}\sum_{j=0}^k x_j$ and $\bar{g}_{k+1} = \frac{1}{\gamma}\bar{x}_k - B\bar{x}_k - \frac{1}{\gamma}\sum_i \omega_i \bar{u}_{i,k+1}$.

**Proposition 2** *We have $g_{k+1} \in \sum_i A_i u_{i,k+1}$. Moreover,*

(i) *If $0 < \inf_{k\in\mathbb{N}} \lambda_k \le \sup_{k\in\mathbb{N}} \lambda_k < \frac{4\beta-\gamma}{2\beta}$, $\big((k+1)\|b_k\|\big)_{k\in\mathbb{N}} \in \ell^1_+$ and $\forall i \in [\![1,n]\!]$, $\big((k+1)\|a_{i,k}\|\big)_{k\in\mathbb{N}} \in \ell^1_+$, then*

$$d\big(0, \textstyle\sum_i A_i u_{i,k+1} + B(\sum_i \omega_i u_{i,k+1})\big) \le \frac{1}{\gamma}\sqrt{\frac{d_0^2 + C_1}{\underline{\tau}(k+1)}}\,.$$

(ii) *If $\underline{\lambda} = \inf_{k\in\mathbb{N}} \lambda_k > 0$, $(\lambda_k\|b_k\|)_{k\in\mathbb{N}} \in \ell^1_+$ and $(\lambda_k\|a_{i,k}\|)_{k\in\mathbb{N}} \in \ell^1_+$, $\forall i \in [\![1,n]\!]$, then*

$$\|\bar{g}_{k+1} + B(\textstyle\sum_i \omega_i \bar{u}_{i,k+1})\| \le \frac{2(d_0 + C_2)}{\gamma\underline{\lambda}(k+1)}.$$

*$C_1, C_2$ are the constants in Theorem 1 and 2 respectively.*

*Proof* (i) Embarking from the definition of $u_{i,k+1}$, we have

$$u_{i,k+1} = J_{\frac{\gamma}{\omega_i}A_i}(2x_k - z_{i,k} - \gamma B x_k)$$
$$\Leftrightarrow \tfrac{\omega_i}{\gamma}\big(2x_k - z_{i,k} - \gamma B x_k - u_{i,k+1} + \gamma B(\textstyle\sum_i \omega_i u_{i,k+1})\big) \in A_i u_{i,k+1} + \omega_i B(\textstyle\sum_i \omega_i u_{i,k+1}).$$

Then sum up over $i$ to get

$$\tfrac{1}{\gamma}x_k - Bx_k - \tfrac{1}{\gamma}\textstyle\sum_i \omega_i u_{i,k+1} + B(\sum_i \omega_i u_{i,k+1}) \in \sum_i A_i u_{i,k+1} + B(\sum_i \omega_i u_{i,k+1}).$$

Now

$$\begin{aligned}\textstyle\sum_{k\in\mathbb{N}}(k+1)\|\boldsymbol{\varepsilon}_k\| &\le \textstyle\sum_{k\in\mathbb{N}}(k+1)\|\boldsymbol{\varepsilon}_{1,k}\| + \sum_{k\in\mathbb{N}}(k+1)\|\boldsymbol{\varepsilon}_{2,k}\| \\ &\le \textstyle\sum_{k\in\mathbb{N}}(k+1)\|b_k\| + \sum_{i=1}^n w_i \sum_{k\in\mathbb{N}}(k+1)\|a_{i,k}\| < +\infty,\end{aligned}$$

and thus Theorem 1 applies. Combining this with the fact that $\mathrm{Id} - \gamma B$ is non-expansive ((ii) of Lemma 4), we obtain

$$\begin{aligned}&d\big(0, \textstyle\sum_i A_i u_{i,k+1} + B(\sum_i \omega_i u_{i,k+1})\big) \\ &\le \|g_{k+1} + \gamma B(\textstyle\sum_i \omega_i u_{i,k+1})\| = \frac{1}{\gamma}\|(\mathrm{Id} - \gamma B)x_k - (\mathrm{Id} - \gamma B)(\sum_i \omega_i u_{i,k+1})\| \\ &\le \frac{1}{\gamma}\|x_k - \textstyle\sum_i \omega_i u_{i,k+1}\| \le \frac{1}{\gamma}\|\boldsymbol{e}_k\| \le \frac{1}{\gamma}\sqrt{\frac{d_0^2 + C_1}{\underline{\tau}(k+1)}}.\end{aligned}$$

(ii) Again, one can check that indeed $(\lambda_k\|\boldsymbol{\varepsilon}_k\|) \in \ell_1^+$. Thus, owing to Theorem 2 and non-expansiveness of $\mathrm{Id}-\gamma B$, we have

$$\begin{aligned}&\|\bar{g}_{k+1} + B(\textstyle\sum_i \omega_i \bar{u}_{i,k+1})\| = \frac{1}{\gamma}\|(\mathrm{Id} - \gamma B)\bar{x}_k - (\mathrm{Id} - \gamma B)(\sum_i \omega_i \bar{u}_{i,k+1})\| \\ &\le \frac{1}{\gamma}\|\bar{x}_k - \textstyle\sum_i \omega_i \bar{u}_{i,k+1}\| \le \frac{1}{\gamma(k+1)}\|\sum_{j=0}^k (\boldsymbol{z}_{j+1} - \boldsymbol{z}_j - \lambda_j \boldsymbol{\varepsilon}_j)/\lambda_j\| \\ &\le \frac{1}{\gamma\underline{\lambda}(k+1)}\big(\|\boldsymbol{z}_0 - \boldsymbol{z}_{k+1}\| + \textstyle\sum_{j=0}^k \lambda_j \|\boldsymbol{\varepsilon}_j\|\big) \le \frac{2(d_0 + C_2)}{\gamma\underline{\lambda}(k+1)}.\end{aligned}$$

*Remark 6*

- Proposition 2 indicates that GFB provides an $\epsilon$-accurate solution in at most $O(1/\epsilon)$ iterations for the criterion $d\big(0, \sum_i A_i u_{i,k} + B(\sum_i \omega_i u_{i,k})\big)^2 \le \epsilon$, which also recovers the sub-differential stopping criterion of the FBS method when $n = 1$. This can then be viewed as a generalization of the best-known complexity bounds of the gradient descent method [51].
- In [11], the Forward–Douglas–Rachford splitting (FDRS) was proposed to solve (1) with $A_1 = A$ and $A_2 = N_{\mathcal{V}}$, $\mathcal{V}$ is a closed vector subspace of $\mathcal{K}$. FDRS has a very similar structure to GFB. Therefore, FDRS also obeys the iteration-complexity bounds established in Section 3.
- Using the transportation formula for cocoercive operators in [63, Lemma 2.8], one can write the *non-relaxed* and *exact* GFB method into the HPE framework introduced there, and derive the iteration-complexity bounds [48] as discussed in the related work. However, we would like to point out that, to establish the iteration-complexity bounds for GFB under the HPE framework, neither relaxation nor errors are handled, and the iteration-complexity bound would be non-uniform [48][4].

### 5.2 Douglas–Rachford splitting and ADMM

As mentioned above, the DRS method [43] can be applied to solve (1) when $B = 0,\ n = 2$ and $L_i = \mathrm{Id},\ i = 1, 2$. Let $\gamma > 0$, $\lambda_k \in ]0, 2]$, $(\varepsilon_{1,k})_{k\in\mathbb{N}}$ and $(\varepsilon_{2,k})_{k\in\mathbb{N}}$ are summable error sequences. The DRS iteration reads

$$\left\lfloor \begin{aligned} u_{k+1} &= J_{\gamma A_1}(2x_k - z_k), \\ z_{k+1} &= z_k + \lambda_k(u_{k+1} + \varepsilon_{1,k} - x_k), \\ v_{k+1} &= J_{\gamma A_2}(z_{k+1}), \\ x_{k+1} &= v_{k+1} + \varepsilon_{2,k+1}. \end{aligned}\right. \tag{19}$$

[4] Let $(x_k, v_k) \in \mathrm{gra}A$ be a sequence generated by an iterative method for solving the monotone inclusion problem $0 \in Ax$, then the non-uniform iteration-complexity bound means that for every $k \in \mathbb{N}$, there exists a $j \le k$ such that $\|v_j\| = O(1/\sqrt{k})$.

When $B = 0,\ n = 2$, $L_2 = \mathrm{Id}$ and $L_1 : \mathcal{H} \to \mathcal{G}$ is some injective linear operator, then (1) can also be solved using ADMM which is DRS method applied to the dual of (1) [27]. In the following, we discuss only DRS.

DRS (19) takes exactly the form (4) with

$$T = \tfrac{1}{2}(R_{\gamma A_1} R_{\gamma A_2} + \mathrm{Id}) \ \text{ and } \ \varepsilon_k = \big(\tfrac{1}{2}(R_{\gamma A_1}(R_{\gamma A_2} z_k + 2\varepsilon_{2,k}) + z_k) - T z_k\big) + \varepsilon_{1,k},$$

see e.g. [18]. Moreover $T \in \mathcal{A}(\frac{1}{2})$ and $(\varepsilon_k)_{k\in\mathbb{N}} \in \ell_+^1$ owing to the summability of $(\varepsilon_{1,k})_{k\in\mathbb{N}}$, $(\varepsilon_{2,k})_{k\in\mathbb{N}}$. Therefore, the DRS iterates obey the iteration-complexity bounds in Section 3.

Next, we turn to the corresponding monotone inclusion problem (1), and develop a criterion similar to Proposition 2. Define $g_{k+1} = \frac{1}{\gamma}(2x_k - z_k - u_{k+1} + z_{k+1} - v_{k+1})$.

**Proposition 3** *We have $g_{k+1} \in A_1 u_{k+1} + A_2 v_{k+1}$. Moreover, if $0 < \inf_{k\in\mathbb{N}} \lambda_k \le \sup_{k\in\mathbb{N}} \lambda_k < 2$, $(\varepsilon_{1,k})_{k\in\mathbb{N}}$ and $(\varepsilon_{2,k})_{k\in\mathbb{N}} \in \ell_+^1$, then*

$$d\big(0, A_1 u_{k+1} + A_2 v_{k+1}\big) \le \frac{1+\lambda_k}{\gamma}\sqrt{\frac{d_0^2 + C_1}{\underline{\tau}(k+1)}} + c_k, \tag{20}$$

*where $c_k = \frac{1}{\gamma}\big((2+\lambda_k)\|\varepsilon_{2,k}\| + \|\varepsilon_{1,k}\|\big)$, $d_0$ and $C_1$ are those defined in Theorem 1.*

*Proof* From (19), we have $2x_k - z_k - u_{k+1} \in \gamma A_1 u_{k+1}, z_{k+1} - v_{k+1} \in \gamma A_2 v_{k+1}$, whose sum leads to

$$g_{k+1} = \tfrac{1}{\gamma}\big(2x_k - z_k - u_{k+1} + z_{k+1} - v_{k+1}\big) \in A_1 u_{k+1} + A_2 v_{k+1},$$

where $\varepsilon_{1,2,k} = \varepsilon_{1,k} + \varepsilon_{2,k}$. Note that $v_k - z_k = (\mathrm{Id} - J_{\gamma A_2}) z_k$ and $(\mathrm{Id} - J_{\gamma A_2})$ is firmly non-expansive (Lemma 1), whence we get

$$d\big(0, A_1 u_{k+1} + A_2 v_{k+1}\big) \le \tfrac{1+\lambda_k}{\gamma}\|e_k\| + \tfrac{1}{\gamma}\big(\lambda_k\|\varepsilon_k\| + \|\varepsilon_{1,2,k} - \varepsilon_k\|\big). \tag{21}$$

For the error $\|\varepsilon_k\|$, we have

$$\|\varepsilon_k\| \le \tfrac{1}{2}\|R_{\gamma A_1}(R_{\gamma A_2} z_k + 2\varepsilon_{2,k}) - R_{\gamma A_1} R_{\gamma A_2} z_k\| + \|\varepsilon_{1,k}\| \le \|\varepsilon_{1,k}\| + \|\varepsilon_{2,k}\|, \tag{22}$$

and similarly we have

$$\|\varepsilon_k - \varepsilon_{1,2,k}\| \le 2\|\varepsilon_{2,k}\|. \tag{23}$$

Combining together (21), (22) and (23) concludes the proof.

*Remark 7*

(1) The summability assumption of $(\varepsilon_{1,k})_{k\in\mathbb{N}}$ and $(\varepsilon_{1,k})_{k\in\mathbb{N}}$ implies $(c_k)_{k\in\mathbb{N}}$ is summable too, hence decays faster than $1/k$, which means the right hand side of (20) is dominated by the first term.

However, to ensure the convergence of the DRS method, one only needs $\big(\lambda_k\|\varepsilon_k\|\big)_{k\in\mathbb{N}}$ to be summable, which does not necessary mean that $(\varepsilon_{1,k})_{k\in\mathbb{N}}$ and $(\varepsilon_{1,k})_{k\in\mathbb{N}}$ should be summable. In [17, Remark 5.7], an example is provided where $\sum_{k\in\mathbb{N}}\|\varepsilon_k\|$ may diverge while $(\lambda_k\|\varepsilon_k\|)_{k\in\mathbb{N}}$ is summable. Suppose $\|\varepsilon_k\| \le (1+\sqrt{1-1/k})/k^q, q \in ]0,1]$, and $\lambda_k = (1-\sqrt{1-1/k})/2$, then it can be verified that $\sum_{k\in\mathbb{N}}\|\varepsilon_k\|$ diverges but $\sum_{k\in\mathbb{N}}\lambda_k\|\varepsilon_k\| < +\infty$ and $\sum_{k\in\mathbb{N}}\lambda_k(2-\lambda_k) = +\infty$.

(2) The obtained iteration-complexity results can be readily adapted to the ADMM method, by exploiting the fact that ADMM is nothing but DRS applied to the Fenchel dual. In particular, one can then show that the pointwise iteration-complexity of ADMM is indeed $O(1/\sqrt{k})$. A similar result in the exact case is presented in [32] under a different metric.

### 5.3 Primal–Dual splitting

In [65], a more general monotone problem is considered. Let $C : \mathcal{H} \to 2^{\mathcal{H}}$ be a maximal monotone operator, $B : \mathcal{H} \to \mathcal{H}$ is $\mu$-cocoercive for some $\mu \in ]0, +\infty[$. $n$ is a strictly positive integer, let $(\omega_i)_i \in ]0,1[^n$ such that $\sum_i \omega_i = 1$. For every $i \in [\![1, n]\!]$, let $\mathcal{G}_i$ be a real Hilbert space, $r_i \in \mathcal{G}_i$, $A_i : \mathcal{G}_i \to 2^{\mathcal{G}_i}$ is maximal monotone, $D_i : \mathcal{G}_i \to 2^{\mathcal{G}_i}$ is maximal monotone and $\nu_i$-strongly monotone for $\nu_i \in ]0, +\infty[$ (thus $D_i^{-1}$ is single-valued), $A_i \Box D_i = (A_i^{-1} + D_i^{-1})^{-1}$ is the parallel sum of $A_i$ and $D_i$. Finally, let $L_i : \mathcal{H} \to \mathcal{G}_i$ be a non-zero bounded linear operator. Now, consider the following monotone inclusion problem,

$$\text{Find } x \in \mathcal{H} \ \text{ such that } \ 0 \in Cx + Bx + \textstyle\sum_{i=1}^n \omega_i L_i^*\big((A_i \Box D_i)(L_i x - r_i)\big), \tag{24}$$

and the corresponding dual problem,

$$\text{Find } (v_1 \in \mathcal{G}_1, ..., v_n \in \mathcal{G}_n) \ \text{ such that } \ (\exists x \in \mathcal{H}) \begin{cases} 0 \in Cx + Bx + \sum_{i=1}^{n} \omega_i L_i^* v_i, \\ 0 \in (A_i \square D_i)(L_i x - r_i) - v_i, \end{cases} \tag{25}$$

denote by $\mathcal{P}$ and $\mathcal{D}$ the solution sets of (24) and (25) respectively.

Let $\tau, (\sigma_i)_i > 0$ such that $\eta = \min\{\frac{1}{\tau}, \frac{1}{\sigma_1}, ..., \frac{1}{\sigma_n}\}(1 - \sqrt{\tau \sum_i \sigma_i \omega_i \|L_i\|^2}), 2\eta\beta > 1$ where $\beta = \min\{\mu, \nu_1, ..., \nu_n\}$, and $\lambda_k \in ]0, \frac{4\eta\beta - 1}{2\eta\beta}]$. Let $(\varepsilon_{1,k}, \varepsilon_{2,k})_{k \in \mathbb{N}}$ be absolutely summable sequences in $\mathcal{H}$, $(\varepsilon_{3,i,k}, \varepsilon_{4,i,k})_{k \in \mathbb{N}}$ be absolutely summable sequences in $\mathcal{G}_i$ for $i \in [\![1, n]\!]$. The primal-dual splitting algorithm of [65] reads,

$$\left\lfloor \begin{array}{l} p_{k+1} = J_{\tau C}\big(x_k - \tau(\sum_i \omega_i L_i^* v_{i,k} + Bx_k + \varepsilon_{1,k})\big) + \varepsilon_{2,k}, \\ y_{k+1} = 2p_{k+1} - x_k, \\ x_{k+1} = x_k + \lambda_k(p_{k+1} - x_k), \\ \text{For } i = 1, \cdots, n \\ \left\lfloor \begin{array}{l} q_{i,k+1} = J_{\sigma_i A_i^{-1}}\big(v_{i,k} + \sigma_i(L_i y_{i,k+1} - D_i^{-1} v_{i,k} - \varepsilon_{3,i,k} - r_i)\big) + \varepsilon_{4,i,k}, \\ v_{i,k+1} = v_{i,k} + \lambda_k(q_{i,k+1} - v_{i,k}). \end{array} \right. \end{array} \right. \tag{26}$$

For the case of convex optimization, with $n = 1, r = 0, D = 0$, no errors and $\lambda_k \equiv 1$, the algorithm reduces to that in [14].

In this part, we recall briefly the fixed-point iteration corresponding to (26) whose detailed derivation can be found in [65, Section 3]. Define the product space $\boldsymbol{\mathcal{G}} = \mathcal{G}_i \times \cdots \times \mathcal{G}_n$ endowed with the scalar product $\langle \boldsymbol{v}_1, \boldsymbol{v}_2 \rangle_{\boldsymbol{\mathcal{G}}} = \sum_{i=1}^{n} \omega_i \langle v_{1,i}, v_{2,i} \rangle_{\mathcal{G}_i}$ and associated norm $\|\cdot\|_{\boldsymbol{\mathcal{G}}}$. Let $\boldsymbol{\mathcal{K}} = \mathcal{H} \oplus \boldsymbol{\mathcal{G}}$ be the Hilbert direct sum with the scalar product $\langle (x_1, \boldsymbol{v}_1), (x_2, \boldsymbol{v}_2) \rangle_{\boldsymbol{\mathcal{K}}} = \langle x_1, x_2 \rangle + \langle \boldsymbol{v}_1, \boldsymbol{v}_2 \rangle_{\boldsymbol{\mathcal{G}}}$ and norm $\|\cdot\|_{\boldsymbol{\mathcal{K}}}$.

Define the following operators on $\boldsymbol{\mathcal{K}}$,

$$\begin{aligned} \boldsymbol{C} &: \boldsymbol{\mathcal{K}} \to 2^{\boldsymbol{\mathcal{K}}}, (x, \boldsymbol{v}) \mapsto (Cx) \times \big(r_1 + A_1^{-1} v_1\big) \times \cdots \times \big(r_n + A_n^{-1} v_n\big), \\ \boldsymbol{D} &: \boldsymbol{\mathcal{K}} \to \boldsymbol{\mathcal{K}}, (x, \boldsymbol{v}) \mapsto \big(\textstyle\sum_i \omega_i L_i^* v_i, -L_1 x, \cdots, -L_n x\big), \\ \boldsymbol{E} &: \boldsymbol{\mathcal{K}} \to \boldsymbol{\mathcal{K}}, (x, \boldsymbol{v}) \mapsto \big(Bx, D_1^{-1} v_1, \cdots, D_n^{-1} v_n\big), \\ \boldsymbol{F} &: \boldsymbol{\mathcal{K}} \to \boldsymbol{\mathcal{K}}, (x, \boldsymbol{v}) \mapsto \big(\tfrac{1}{\tau} x - \textstyle\sum_i \omega_i L_i^* v_i, \tfrac{1}{\sigma_1} v_1 - L_1 x, \cdots, \tfrac{1}{\sigma_n} v_n - L_n x\big), \end{aligned}$$

$\boldsymbol{C}$ and $\boldsymbol{D}$ are maximal monotone, $\boldsymbol{E}$ is $\beta$–cocoercive, and $\boldsymbol{F}$ is self–adjoint and $\eta$–strongly positive. Then the fixed point equation of (26) is [65]

$$\boldsymbol{z}_{k+1} = \boldsymbol{z}_k + \lambda_k\big(J_{\boldsymbol{A}}(\boldsymbol{z}_k - \boldsymbol{B}\boldsymbol{z}_k - \boldsymbol{\varepsilon}_{2,k}) + \boldsymbol{\varepsilon}_{1,k} - \boldsymbol{z}_k\big), \tag{27}$$

where $\boldsymbol{A} = \boldsymbol{F}^{-1}(\boldsymbol{C} + \boldsymbol{D}), \boldsymbol{B} = \boldsymbol{F}^{-1}\boldsymbol{E}, \boldsymbol{z}_k = (x_k, v_{1,k}, \cdots, v_{n,k})$, and the errors $\boldsymbol{\varepsilon}_{1,k} = (\varepsilon_{2,k}, \varepsilon_{4,1,k}, \cdots, \varepsilon_{4,n,k})$, $\boldsymbol{\varepsilon}_{3,k} = (\varepsilon_{1,k}, \varepsilon_{3,1,k}, \cdots, \varepsilon_{3,n,k}), \boldsymbol{\varepsilon}_{4,k} = \big(\frac{1}{\tau}\varepsilon_{2,k}, \frac{1}{\sigma_1}\varepsilon_{4,1,k}, \cdots, \frac{1}{\sigma_n}\varepsilon_{4,n,k}\big), \boldsymbol{\varepsilon}_{2,k} = \boldsymbol{F}^{-1}\big((\boldsymbol{D} + \boldsymbol{E})\boldsymbol{\varepsilon}_{1,k} + \boldsymbol{\varepsilon}_{3,k} - \boldsymbol{\varepsilon}_{4,k}\big)$.

Iteration (27) has the structure of FBS, and by Remark 1, the fixed-point operator $\boldsymbol{T} = J_{\boldsymbol{A}} \circ (\mathbf{Id} - \boldsymbol{B}) \in \mathcal{A}(\frac{2\eta\beta}{4\eta\beta - 1})$. Thus, (27) is a special instance of (4). In addition, since $\varepsilon_{1,k}$ and $\varepsilon_{2,k}$ (resp. $\varepsilon_{3,i,k}$ and $\varepsilon_{4,i,k}$) are summable in $\mathcal{H}$ (resp. in $\mathcal{G}_i$), so are $\boldsymbol{\varepsilon}_{1,k}$ and $\boldsymbol{\varepsilon}_{2,k}$ in $\boldsymbol{\mathcal{K}}$. Therefore, the PDS iterates obey the iteration-complexity bounds in Section 3.

Observing that $\mathrm{fix}\boldsymbol{T} = \mathrm{zer}(\boldsymbol{A} + \boldsymbol{B}) = \mathrm{zer}(\boldsymbol{C} + \boldsymbol{D} + \boldsymbol{E})$, we also have the following bounds. We will denote $\delta = \max\{\frac{1}{\tau}, \frac{1}{\sigma_1}, \ldots, \frac{1}{\sigma_n}\}$, $\underline{\tau} = \inf_{k \in \mathbb{N}} \lambda_k(\frac{4\eta\beta - 1}{2\eta\beta} - \lambda_k)$, and $\boldsymbol{w}_{k+1} = (\boldsymbol{z}_{k+1} - (1 - \lambda_k)\boldsymbol{z}_k)/\lambda_k - \boldsymbol{\varepsilon}_k$, where $\boldsymbol{\varepsilon}_k = (J_{\boldsymbol{A}}(\boldsymbol{z}_k - \boldsymbol{B}\boldsymbol{z}_k - \boldsymbol{\varepsilon}_{2,k}) + \boldsymbol{\varepsilon}_{1,k}) - J_{\boldsymbol{A}}(\boldsymbol{z}_k - \boldsymbol{B}\boldsymbol{z}_k)$.

**Proposition 4** *Suppose* $0 < \inf_{k \in \mathbb{N}} \lambda_k \leq \sup_{k \in \mathbb{N}} \lambda_k < \frac{2\eta\beta}{4\eta\beta - 1}$, $\big((k+1)\|\varepsilon_{j,k}\|\big)_{k \in \mathbb{N}} \in \ell_+^1$, $j = 1, 2$, *and* $\big((k+1)\|\varepsilon_{j,i,k}\|\big)_{k \in \mathbb{N}} \in \ell_+^1$, $j = 3, 4$. *Then*

$$d\big(0, (\boldsymbol{A} + \boldsymbol{B})\boldsymbol{w}_{k+1}\big) \leq \frac{2\delta}{\eta}\sqrt{\frac{d_0^2 + C_1}{\underline{\tau}(k+1)}},$$

*where* $C_1 < +\infty$ *and* $d_0 = \inf_{\boldsymbol{z} \in \mathrm{fix}\boldsymbol{T}} \|\boldsymbol{z}_0 - \boldsymbol{z}\|_{\boldsymbol{\mathcal{K}}}$.

The proof is similar to that of Proposition 2 and will not be included here. Again, Proposition 4 can serve as a stopping criterion for the primal-dual monotone inclusion (24)-(25). An ergodic bound can also be derived. However, for the sake of brevity, we do not pursue this further.

## 6 Non-stationary Krasnosel'skiĭ-Mann iteration

### 6.1 General convergence analysis

The fixed-point iteration discussed in Section 3 is stationary, namely, $T$ of (4) is fixed along the iterations. In this section, we study the non-stationary version of it, and show that, the non-stationary case can be seen as a perturbation of the stationary one. Moreover the iterates are convergent if the extra perturbation error is absolutely summable.

Let $T_\Gamma : \mathcal{H} \to \mathcal{H}$ be a non-expansive operator depending on a parameter $\Gamma$. Let $\lambda_k \in ]0,1[$. Then the non-stationary fixed-point iteration is defined by

$$z_{k+1} = z_k + \lambda_k(T_{\Gamma_k} z_k + \varepsilon_k - z_k) = T_{\Gamma_k,\lambda_k} z_k + \lambda_k \varepsilon_k, \tag{28}$$

with $T_{\Gamma_k,\lambda_k} = \lambda_k T_{\Gamma_k} + (1-\lambda_k)\mathrm{Id}$. If we define $\varepsilon_{\Gamma_k} = (T_{\Gamma_k} - T_\Gamma)z_k$, $\pi_k = \varepsilon_{\Gamma_k} + \varepsilon_k$, then (28) can be rewritten as

$$z_{k+1} = \big(\lambda_k T_\Gamma + (1-\lambda_k)\mathrm{Id}\big) z_k + \lambda_k \pi_k = T_{\Gamma,\lambda_k} z_k + \lambda_k \pi_k, \tag{29}$$

where $T_{\Gamma,\lambda_k} = \lambda_k T_\Gamma + (1-\lambda_k)\mathrm{Id}$. The corresponding vector $e_k$ of (28) becomes

$$e_k = (z_k - z_{k+1})/\lambda_k + \pi_k.$$

Comparing (29) to Definition 5, a new error sequence $\pi_k$ is introduced. To obtain convergence of the non-stationary iteration, we adapt arguments from [36] (Banach spaces endowed with an appropriate compatible topology) and [37] (real Hilbert spaces). For convenience, we recall that $\underline{\tau} = \inf_{k\in\mathbb{N}} \lambda_k(1-\lambda_k)$.

**Theorem 4 (Convergence of (28))** *Assume the following holds:*

(A.1) $\mathrm{fix}T_\Gamma \neq \emptyset$;
(A.2) $\forall k \in \mathbb{N},\ T_{\Gamma_k,\lambda_k}$ *is* $(1+\beta_k)$*-Lipschitz with* $\beta_k \geq 0$, *and* $(\beta_k)_{k\in\mathbb{N}} \in \ell^1_+$;
(A.3) $\lambda_k \in ]0,1[$ *such that* $\underline{\tau} > 0$;
(A.4) $(\lambda_k \|\varepsilon_k\|) \in \ell^1_+$;
(A.5) $\forall \rho \in [0,+\infty[$, *the sequence* $(\lambda_k \Delta_{k,\rho})_{k\in\mathbb{N}} \in \ell^+_1$, *where*

$$\Delta_{k,\rho} = \sup_{\|z\|\leq\rho} \|T_{\Gamma_k} z - T_\Gamma z\|. \tag{30}$$

*Then* $(e_k)_{k\in\mathbb{N}}$ *converges strongly to* 0, *and* $(z_k)_{k\in\mathbb{N}}$ *converges weakly to a point* $z^\star \in \mathrm{fix}T_\Gamma$.

*Proof* For the sequence $(z_k)_{k\in\mathbb{N}}$ generated by (28), and $z^\star \in \mathrm{fix}T_\Gamma$, we have

$$\begin{aligned}\|z_{k+1} - z^\star\| &= \|T_{\Gamma_k,\lambda_k} z_k + \lambda_k \varepsilon_k - T_{\Gamma,\lambda_k} z^\star\| \\ &\leq \|T_{\Gamma_k,\lambda_k} z_k - T_{\Gamma_k,\lambda_k} z^\star\| + \|T_{\Gamma_k,\lambda_k} z^\star - T_{\Gamma,\lambda_k} z^\star\| + \lambda_k\|\varepsilon_k\| \\ &\leq (1+\beta_k)\|z_k - z^\star\| + \lambda_k \Delta_{k,\|z^\star\|} + \lambda_k\|\varepsilon_k\|.\end{aligned}$$

As $(\beta_k)_{k\in\mathbb{N}}$, $(\lambda_k\|\varepsilon_k\|)_{k\in\mathbb{N}}$ and $(\lambda_k \Delta_{k,\|z^\star\|})_{k\in\mathbb{N}}$ are summable by assumptions (A.2), (A.4) and (A.5), it follows from [57, Lemma 2.2.2] that the sequence $(\|z_k - z^\star\|)_{k\in\mathbb{N}}$ converges, hence bounded. Therefore, $z_k$ is bounded in norm by some $\rho \in [0,+\infty[$. This implies that

$$\|z_{k+1} - T_{\Gamma,\lambda_k} z_k\| = \|T_{\Gamma_k,\lambda_k} z_k + \lambda_k\varepsilon_k - T_{\Gamma,\lambda_k} z_k\| \leq \lambda_k \left(\Delta_{k,\rho} + \|\varepsilon_k\|\right).$$

In other words, the (inexact) non-stationary iteration (28) can be seen as a perturbed version of the (inexact) stationary one with an extra-error term which is summable owing to (A.5). The rest of the proof follows by applying [36, Proposition 2.1 and Remark 2.2] (see also [10, Remark 14]) using (A.1) and (A.3).

If $T_{\Gamma_k}$ were non-expansive, then obviously (A.2) is in force, and in turn, Theorem 4 holds. In the specific scenario where $\forall k \in \mathbb{N},\ T_{\Gamma_k}$ is $\alpha_k$-averaged, $\alpha_k \in ]0,1]$, we can further refine the choice of $\lambda_k$. Here we take a different route from the one in [18]. By definition, $\forall k \in \mathbb{N}$, there exists a non-expansive operator $R_{\Gamma_k} : \mathcal{H} \to \mathcal{H}$ such that $T_{\Gamma_k} = \alpha_k R_{\Gamma_k} + (1-\alpha_k)\mathrm{Id}$. Let $R_\Gamma$ be a non-expansive operator, and $\lambda'_k = \alpha_k\lambda_k$. We have the following corollary.

**Corollary 5** *Assume that*

(A′.1) $\text{fix} R_\Gamma \neq \emptyset$;
(A′.2) $\forall k \in \mathbb{N},\ T_{\Gamma_k}$ *is* $\alpha_k$*-averaged,* $\alpha_k \in ]0, 1]$;
(A′.3) $\lambda_k \in ]0, \frac{1}{\alpha_k}[$ *such that* $\inf_{k\in\mathbb{N}} \lambda'_k(1 - \lambda'_k) > 0$;
(A′.4) $(\lambda'_k \|\varepsilon_k\|) \in \ell^1_+$;
(A′.5) $\forall \rho \in [0, +\infty[,\ (\lambda'_k \Delta_{k,\rho})_{k\in\mathbb{N}} \in \ell^+_1$, *where*

$$\Delta_{k,\rho} = \sup_{\|z\|\leq\rho} \|R_{\Gamma_k} z - R_\Gamma z\|.$$

*Then* $(e_k)_{k\in\mathbb{N}}$ *converges strongly to* $0$, *and* $(z_k)_{k\in\mathbb{N}}$ *converges weakly to a fixed point* $z^\star \in \text{fix} T_\Gamma$.

Assumptions (A.4)-(A.5) of Theorem 4 imply that

$$\textstyle\sum_{k\in\mathbb{N}} \lambda_k \|\pi_k\| \leq \sum_{k\in\mathbb{N}} \lambda_k (\|\varepsilon_{\Gamma_k}\| + \|\varepsilon_k\|) < +\infty.$$

Therefore, if we can further impose a stronger summability assumption on $(\pi_k)_{k\in\mathbb{N}}$ as in Theorem 1, then we can obtain the iteration-complexity bounds for the non-stationary iteration (28) as well. This is stated in the following result. Recall $\Lambda_k = \sum_{j=0}^k \lambda_j$, $\bar{e}_k = \frac{1}{\Lambda_k} \sum_{j=0}^k \lambda_j e_j$, and $d_0 = \inf_{z\in\text{fix}T_\Gamma} \|z_0 - z\|$.

**Proposition 5** *Assume that* (A.1) *holds, and that* $\forall k \in \mathbb{N},\ T_{\Gamma_k}$ *is non-expansive.*

(i) *Suppose that* (A.3) *is verified,* $\big((k + 1)\|\varepsilon_k\|\big)_{k\in\mathbb{N}} \in \ell^+_1$ *and* $\big((k + 1)\Delta_{k,\rho}\big)_{k\in\mathbb{N}} \in \ell^+_1$, $\forall \rho \in [0, +\infty[$, *where* $\Delta_{k,\rho}$ *is given by* (30). *Then,*

$$\|e_k\| \leq \sqrt{\frac{d_0^2 + C_1}{\underline{\tau}(k + 1)}}$$

*where* $C_1$ *is a bounded constant (see Theorem 1).*

(ii) *Suppose that* $\lambda_k \in ]0, 1]$ *such that* (A.4)-(A.5) *are verified. Then*

$$\|\bar{e}_k\| \leq \frac{2(d_0 + C_2)}{\Lambda_k}$$

*where* $C_2$ *is a bounded constant (see Theorem 2). In particular, if* $\inf_{k\in\mathbb{N}} \lambda_k > 0$, *we get* $O(1/k)$ *ergodic convergence rate.*

*Proof* (i) By assumption, we have $\sum_{k\in\mathbb{N}} (k+1)\|\pi_k\| < +\infty$. All assumptions of Theorem 1 are then fulfilled and the result follows. (ii) Similarly all required assumptions to apply Theorem 2 are in force.

*Remark 8* If metric sub-regularity assumption is imposed on $\text{Id} - T_\Gamma$, a result similar to Theorem 3 can be stated. But now, we have $c_k = \nu_1 \lambda_k (\|\varepsilon_k\| + \Delta_{k,\rho})$. So the actual local convergence behavior depends also on the additional perturbation error brought by non-stationarity as captured by $\Delta_{k,\rho}$, even in the exact case. Thus, similarly to Remark 4, if $\|\varepsilon_k\|$ and $\Delta_{k,\rho}$ converge linearly, then so does $d_k = d(z_k, \text{fix} T_\Gamma)$ locally. If the non-stationary error $\Delta_{k,\rho}$ decays sub-linearly, then it dominates.

### 6.2 Application to non-stationary GFB

As discussed in Section 5.1, the fixed-point operator $\boldsymbol{T}_\gamma$ of GFB depends on a parameter $\gamma$. Now let $\gamma$ vary along the iteration, and denote the following operators

$$\boldsymbol{T}_{1,\gamma_k} = \tfrac{1}{2}(R_{\gamma_k \boldsymbol{A}} R_{\boldsymbol{\mathcal{S}}} + \mathbf{Id}),\ \boldsymbol{T}_{2,\gamma_k} = \mathbf{Id} - \gamma_k \boldsymbol{B}_{\boldsymbol{\mathcal{S}}},\ \text{ and } \boldsymbol{T}_{\gamma_k} = \boldsymbol{T}_{1,\gamma_k} \circ \boldsymbol{T}_{2,\gamma_k}.$$

Recall that for $(\gamma_k)_{k\in\mathbb{N}} \in ]0, 2\beta[$, $\boldsymbol{T}_{\gamma_k}$ is $\alpha_k$-averaged with $\alpha_k = \frac{2\beta}{4\beta-\gamma_k}$, and $\boldsymbol{T}_\gamma$ is $\alpha$-averaged with $\alpha = \frac{2\beta}{4\beta-\gamma}$, for $\gamma \in ]0, 2\beta[$. Moreover, there exists non-expansive operators $\boldsymbol{R}_\gamma$ and $\boldsymbol{R}_{\gamma_k}$ such that $\boldsymbol{T}_\gamma = \alpha \boldsymbol{R}_\gamma + (1 - \alpha)\mathbf{Id}$ and $\boldsymbol{T}_{\gamma_k} = \alpha_k \boldsymbol{R}_{\gamma_k} + (1 - \alpha_k)\mathbf{Id}$.

Let $\boldsymbol{\varepsilon}_{\gamma_k} = (\boldsymbol{T}_{\gamma_k} - \boldsymbol{T}_\gamma)\boldsymbol{z}_k,\ \boldsymbol{\pi}_k = \boldsymbol{\varepsilon}_{\gamma_k} + \boldsymbol{\varepsilon}_k$. The the non-stationary version of (18) is defined by

$$\boldsymbol{z}_{k+1} = \boldsymbol{z}_k + \lambda_k (\boldsymbol{T}_\gamma \boldsymbol{z}_k + \boldsymbol{\pi}_k - \boldsymbol{z}_k). \tag{31}$$

**Theorem 5** *For the non-stationary iteration* (31)*, if the following assumptions hold*

(A″.1) $\text{zer}(B + \sum_i A_i) \neq \emptyset$;

(A$''$.2) $\lambda_k \in ]0, \frac{1}{\alpha_k}[$, *such that* $\inf_{k\in\mathbb{N}} \lambda_k(\frac{1}{\alpha_k} - \lambda_k) > 0$*;*
(A$''$.3) $\big(\lambda_k\|b_k\|\big)_{k\in\mathbb{N}} \in \ell_+^1$ *and* $\big(\lambda_k\|a_{i,k}\|\big)_{k\in\mathbb{N}} \in \ell_+^1$, $\forall i \in [\![1,n]\!]$*;*
(A$''$.4) $(\gamma_k)_{k\in\mathbb{N}} \in ]0, 2\beta[$ *such that* $0 < \underline{\gamma} \le \gamma_k \le \overline{\gamma} < 2\beta$, $\gamma \in [\underline{\gamma}, \overline{\gamma}]$, *and* $(\lambda_k|\gamma_k - \gamma|)_{k\in\mathbb{N}} \in \ell_1^+$*;*

*then, the sequence* $(\boldsymbol{z}_k)_{k\in\mathbb{N}}$ *generated by* (31) *converges weakly to a point in* $\mathrm{fix}\boldsymbol{T}_{\boldsymbol{\gamma}}$.
*If we further assume that* $((k+1)|\gamma_k - \gamma|)_{k\in\mathbb{N}} \in \ell_1^+$, $((k+1)\|b_k\|)_{k\in\mathbb{N}} \in \ell_+^1$, *and* $((k+1)\|a_{i,k}\|)_{k\in\mathbb{N}} \in \ell_+^1$, $\forall i \in [\![1,n]\!]$, *then we obtain the pointwise iteration-complexity bound for the non-stationary version of GFB algorithm as stated in (i) of Proposition 5.*

*Proof* It is sufficient to verify the conditions of Corollary 5 to conclude.

- Assumption (A$'$.1) is fulfilled thanks to (A$''$.1) since $P_{\boldsymbol{\mathcal{S}}}(\mathrm{fix}\boldsymbol{R}_{\boldsymbol{\gamma}}) = P_{\boldsymbol{\mathcal{S}}}(\mathrm{fix}\boldsymbol{T}_{\boldsymbol{\gamma}}) = \mathrm{zer}(B + \sum_i A_i) \ne \emptyset$;
- As $\boldsymbol{T}_{\gamma_k} \in \mathcal{A}(\alpha_k)$, $\boldsymbol{T}_{\gamma_k,\lambda_k} \in \mathcal{A}(\alpha_k\lambda_k)$, and thus assumption (A$'$.2) is in force;
- (A$'$.3) holds thanks to (A$''$.2);
- (A$'$.4) follows from (A$''$.3).
- It remains to check that (A$''$.4) implies (A$'$.5).

By definition of $\boldsymbol{R}_{\boldsymbol{\gamma}}$ and $\boldsymbol{R}_{\gamma_k}$, we have

$$\boldsymbol{R}_{\boldsymbol{\gamma}} = (1 - \tfrac{1}{\alpha})\mathrm{Id} + \tfrac{1}{\alpha}\boldsymbol{T}_{\boldsymbol{\gamma}} \quad \text{and} \quad \boldsymbol{R}_{\gamma_k} = (1 - \tfrac{1}{\alpha_k})\mathrm{Id} + \tfrac{1}{\alpha_k}\boldsymbol{T}_{\gamma_k} \ .$$

It then follows that

$$\begin{aligned}
\|\boldsymbol{R}_{\gamma_k}\boldsymbol{z} - \boldsymbol{R}_{\boldsymbol{\gamma}}\boldsymbol{z}\| &= \|\tfrac{1}{\alpha_k}(\boldsymbol{T}_{\gamma_k} - \mathrm{Id})\boldsymbol{z} - \tfrac{1}{\alpha}(\boldsymbol{T}_{\boldsymbol{\gamma}} - \mathrm{Id})\boldsymbol{z}\| \\
&\le |\tfrac{1}{\alpha_k} - \tfrac{1}{\alpha}|\|(\boldsymbol{T}_{\boldsymbol{\gamma}} - \mathrm{Id})\boldsymbol{z}\| + \tfrac{1}{\alpha_k}\|(\boldsymbol{T}_{\gamma_k} - \mathrm{Id})\boldsymbol{z} - (\boldsymbol{T}_{\boldsymbol{\gamma}} - \mathrm{Id})\boldsymbol{z}\| \\
&\le \frac{|\gamma_k - \gamma|}{2\beta}(2\rho + \|\boldsymbol{T}_{\boldsymbol{\gamma}}0\|) + \tfrac{1}{\alpha_k}\|\boldsymbol{T}_{\gamma_k}\boldsymbol{z} - \boldsymbol{T}_{\boldsymbol{\gamma}}\boldsymbol{z}\| \ .
\end{aligned} \tag{32}$$

Now, non-expansiveness of $\boldsymbol{T}_{1,\gamma_k}$ yields

$$\begin{aligned}
\|\boldsymbol{T}_{\gamma_k}\boldsymbol{z} - \boldsymbol{T}_{\boldsymbol{\gamma}}\boldsymbol{z}\| &\le \|\boldsymbol{T}_{1,\gamma_k}\boldsymbol{T}_{2,\gamma_k}\boldsymbol{z} - \boldsymbol{T}_{1,\gamma_k}\boldsymbol{T}_{2,\gamma}\boldsymbol{z}\| + \|\boldsymbol{T}_{1,\gamma_k}\boldsymbol{T}_{2,\gamma}\boldsymbol{z} - \boldsymbol{T}_{1,\gamma}\boldsymbol{T}_{2,\gamma}\boldsymbol{z}\| \\
&\le \underbrace{\|\boldsymbol{T}_{2,\gamma_k}\boldsymbol{z} - \boldsymbol{T}_{2,\gamma}\boldsymbol{z}\|}_{\text{Term 1}} + \underbrace{\|\boldsymbol{T}_{1,\gamma_k}\boldsymbol{T}_{2,\gamma}\boldsymbol{z} - \boldsymbol{T}_{1,\gamma}\boldsymbol{T}_{2,\gamma}\boldsymbol{z}\|}_{\text{Term 2}}.
\end{aligned} \tag{33}$$

We first bound the first term in (33),

$$\begin{aligned}
\|\boldsymbol{T}_{2,\gamma_k}\boldsymbol{z} - \boldsymbol{T}_{2,\gamma}\boldsymbol{z}\| &\le |\gamma_k - \gamma|\|\boldsymbol{B}P_{\boldsymbol{\mathcal{S}}}\boldsymbol{z}\| \\
\text{(Triangle inequality and } B \text{ is } \beta^{-1}\text{-Lipschitz)} &\le (\beta^{-1}\rho + \|B(0)\|)|\gamma_k - \gamma|,
\end{aligned} \tag{34}$$

where $B(0)$ is bounded.

Let's now turn to the second term of (33). Denote $\boldsymbol{z}_{\boldsymbol{\mathcal{S}}} = P_{\boldsymbol{\mathcal{S}}}(\boldsymbol{z})$ and $\boldsymbol{z}_{\boldsymbol{\mathcal{S}}^\perp} = \boldsymbol{z} - \boldsymbol{z}_{\boldsymbol{\mathcal{S}}}$, then

$$\boldsymbol{v} = \boldsymbol{T}_{1,\gamma}\boldsymbol{T}_{2,\gamma'}\boldsymbol{z} \Leftrightarrow \boldsymbol{v} = \boldsymbol{z}_{\boldsymbol{\mathcal{S}}^\perp} + J_{\gamma\boldsymbol{A}}(\boldsymbol{z}_{\boldsymbol{\mathcal{S}}} - \boldsymbol{z}_{\boldsymbol{\mathcal{S}}^\perp} - \gamma'\boldsymbol{B}\boldsymbol{z}_{\boldsymbol{\mathcal{S}}}),$$

Let $\boldsymbol{y} = \boldsymbol{z}_{\boldsymbol{\mathcal{S}}} - \boldsymbol{z}_{\boldsymbol{\mathcal{S}}^\perp} - \gamma\boldsymbol{B}\boldsymbol{z}_{\boldsymbol{\mathcal{S}}}$, then we have

$$\boldsymbol{T}_{1,\gamma_k}\boldsymbol{T}_{2,\gamma}\boldsymbol{z} - \boldsymbol{T}_{1,\gamma}\boldsymbol{T}_{2,\gamma}\boldsymbol{z} = J_{\gamma_k\boldsymbol{A}}(\boldsymbol{y}) - J_{\gamma\boldsymbol{A}}(\boldsymbol{y}).$$

and denote $\boldsymbol{u}_k = J_{\gamma_k\boldsymbol{A}}\boldsymbol{y}$, $\boldsymbol{u} = J_{\gamma\boldsymbol{A}}\boldsymbol{y}$. By definition of the resolvent, this is equivalent to

$$\big(\boldsymbol{u}, \frac{\boldsymbol{y} - \boldsymbol{u}}{\gamma}\big) \in \mathrm{gra}(\boldsymbol{A}) \quad \text{and} \quad \big(\boldsymbol{u}_k, \frac{\boldsymbol{y} - \boldsymbol{u}_k}{\gamma_k}\big) \in \mathrm{gra}(\boldsymbol{A}).$$

Since $\boldsymbol{A}$ is monotone, and by assumptions on $\gamma_k$ and $\gamma$, it follows that

$$\langle\gamma_k(\boldsymbol{y} - \boldsymbol{u}) - \gamma(\boldsymbol{y} - \boldsymbol{u}_k), \boldsymbol{u} - \boldsymbol{u}_k\rangle \ge 0 \Leftrightarrow \|\boldsymbol{u} - \boldsymbol{u}_k\|^2 \le \frac{\gamma_k - \gamma}{\gamma}\langle\boldsymbol{y} - \boldsymbol{u}, \boldsymbol{u} - \boldsymbol{u}_k\rangle.$$

Therefore, using Cauchy-Schwartz inequality and the fact that $\mathbf{Id} - J_{\gamma\boldsymbol{A}} \in \mathcal{A}(\frac{1}{2})$ is non-expansive,

$$\|\boldsymbol{u} - \boldsymbol{u}_k\| \le \frac{|\gamma_k - \gamma|}{\underline{\gamma}}\|(\mathbf{Id} - J_{\gamma\boldsymbol{A}})\boldsymbol{y}\| \le \frac{|\gamma_k - \gamma|}{\underline{\gamma}}(\|\boldsymbol{y}\| + \|J_{\gamma\boldsymbol{A}}(\mathbf{0})\|), \tag{35}$$

where $\|J_{\gamma \boldsymbol{A}}(0)\|^2 = \sum_i \omega_i \|J_{\frac{\gamma}{\omega_i}}(0)\|^2$. Using the triangle inequality, the Pythagorean theorem and non-expansiveness of $\beta \boldsymbol{B} P_{\boldsymbol{\mathcal{S}}}$, we obtain

$$\begin{aligned}\|\boldsymbol{y}\| &\le \|\boldsymbol{z}_{\boldsymbol{\mathcal{S}}} - \boldsymbol{z}_{\boldsymbol{\mathcal{S}}^\perp}\| + \gamma \|\boldsymbol{B}\boldsymbol{z}_{\boldsymbol{\mathcal{S}}}\| \le \rho + \gamma \|\boldsymbol{B}\boldsymbol{z}_{\boldsymbol{\mathcal{S}}} - \boldsymbol{B}P_{\boldsymbol{\mathcal{S}}}(0)\| + \gamma \|\boldsymbol{B}P_{\boldsymbol{\mathcal{S}}}(0)\| \\ &\le \rho + \gamma\beta^{-1}\|\boldsymbol{z}\| + \gamma\|B(0)\| \le \rho + \overline{\gamma}\beta^{-1}\rho + \overline{\gamma}\|B(0)\|. \end{aligned} \tag{36}$$

Putting together (32), (34), (35) and (36), we get $\forall \rho \in [0, +\infty[$

$$\sum_{k\in\mathbb{N}} \lambda_k \alpha_k \Delta_{k,\rho} = \sum_{k\in\mathbb{N}} \lambda_k \alpha_k \sup_{\|\boldsymbol{z}\| \le \rho} \|\boldsymbol{R}_{\gamma_k}\boldsymbol{z} - \boldsymbol{R}_{\gamma}\boldsymbol{z}\| \le C \sum_{k\in\mathbb{N}} \lambda_k |\gamma_k - \gamma| < +\infty,$$

where

$$C = \frac{2\rho + \|\boldsymbol{T}_\gamma 0\|}{4\beta - \overline{\gamma}} + \rho\beta^{-1}(1 + \beta/\underline{\gamma} + \overline{\gamma}/\underline{\gamma}) + (1 + \overline{\gamma}/\underline{\gamma})\|B(0)\| + \underline{\gamma}^{-1}\|J_{\gamma \boldsymbol{A}}(\boldsymbol{0})\| < +\infty\ .$$

Consequently, (A.5) is fulfilled.

The last statement of the theorem is a simple application of Theorem 1.

*Remark 9*

(1) For the non-stationary versions of the methods discussed in Section 5, for instance, the DRS and FDRS methods, whose fixed-point operators also depend on $\gamma_k$, Theorem 5 is also applicable. In general, the summability assumption (A$''$.4) of Theorem 5 is hard to remove, except for some special cases. For instance for the FBS method, as stated in [20, Theorem 3.4], where $\gamma_k \in ]0, 2\beta[$, but $\lambda_k \in ]0, 1[$ instead of $]0, \frac{4\beta - \gamma_k}{2\beta}[$ here.

(2) Recall the term $c_k$ in the bound (13) of Theorem 3. Clearly, for the non-stationary GFB iteration, even if the approximation error $\varepsilon_k = 0$, we still have $c_k = \nu_1 \lambda_k \varepsilon_{\gamma_k} \neq 0$. Therefore, under metric sub-regularity of $\mathrm{Id} - \boldsymbol{T}_\gamma$, a bound similar to (13) can be obtained, whose performance will depend on how fast $|\gamma_k - \gamma|$ converges to 0.

## 7 Conclusion

In this paper, we presented global iteration-complexity bounds for the inexact Krasnosel'skiĭ–Mann iteration built from a non-expansive operator, then under metric sub-regularity, we provided a unified quantitative analysis of local linear convergence. Extensions to the non-stationary version of the fixed-point iteration were also studied. The obtained results are applied to several monotone operator splitting algorithms and illustrated through some examples where both global sub-linear and local linear convergence profiles are observed. The local linear convergence rate depends on the sub-regularity modulus of the fixed-point operator, which is not straightforward to compute in general. Moreover, our rate estimates are not necessarily sharp in general. Sharper rates can be obtained if more structure on the problem is available and exploited wisely. These are important aspects that we will investigate in a future work.

## Acknowledgements

This work has been partly supported by the European Research Council (ERC project SIGMA-Vision). J. Fadili is partly supported by Institut Universitaire de France. We would like to thank Yuchao Tang for pointing reference [52] to us.

## References

1. F. J. Aragòn and M. H. Geoffroy. Characterization of metric regularity of subdifferentials. *Journal of Convex Analysis*, 15(2):365–380, 2008.
2. J. B. Baillon. Un théoreme de type ergodique pour les contractions non linéaires dans un espace de Hilbert. *CR Acad. Sci. Paris Sér. AB*, 280:1511–1514, 1975.
3. J. B. Baillon and R. E. Bruck. The rate of asymptotic regularity is $O(1/\sqrt{n})$. *Lecture Notes in Pure and Applied Mathematics*, 178:51–81, 1996.
4. J. B. Baillon and G. Haddad. Quelques propriétés des opérateurs angle-bornés etn-cycliquement monotones. *Israel Journal of Mathematics*, 26(2):137–150, 1977.

5. H. H. Bauschke, J.Y. Bello Cruz, T.A. Nghia, H. M. Phan, and X. Wang. Optimal rates of convergence of matrices with applications. arxiv:1407.0671, 2014.
6. H. H. Bauschke and P. L. Combettes. *Convex analysis and monotone operator theory in Hilbert spaces*. Springer, 2011.
7. H. H. Bauschke, D. R. Luke, H. M. Phan, and X. Wang. Restricted normal cones and the method of alternating projections: theory. *Set-Valued and Variational Analysis*, 21(3):431–473, 2013.
8. H. H. Bauschke, D. Noll, and H. M. Phan. Linear and strong convergence of algorithms involving averaged nonexpansive operators. *Journal of Mathematical Analysis and Applications*, 421(1):1–20, 2015.
9. J. M. Borwein and B. Sims. The Douglas–Rachford algorithm in the absence of convexity. In *Fixed-Point Algorithms for Inverse Problems in Science and Engineering*, pages 93–109. Springer, 2011.
10. H. Brézis and P. L. Lions. Produits infinis de résolvantes. *Israel Journal of Mathematics*, 29(4):329–345, 1978.
11. L. M. Briceno-Arias. Forward–Douglas–Rachford splitting and forward-partial inverse method for solving monotone inclusions. *Optimization*, 64(5):1239–1261, 2015.
12. L. M. Briceno-Arias and P. L. Combettes. A monotone+ skew splitting model for composite monotone inclusions in duality. *SIAM Journal on Optimization*, 21(4):1230–1250, 2011.
13. R. S. Burachik, S. Scheimberg, and B. F. Svaiter. Robustness of the hybrid extragradient proximal-point algorithm. *Journal of Optimization Theory and Applications*, 111(1):117–136, 2001.
14. A. Chambolle and T. Pock. A first–order primal–dual algorithm for convex problems with applications to imaging. *Journal of Mathematical Imaging and Vision*, 40(1):120–145, 2011.
15. G. Chen and M. Teboulle. A proximal-based decomposition method for convex minimization problems. *Mathematical Programming*, 64(1-3):81–101, 1994.
16. P. Chen, J. Huang, and X. Zhang. A primal–dual fixed point algorithm for convex separable minimization with applications to image restoration. *Inverse Problems*, 29(2):025011, 2013.
17. P. L. Combettes. Quasi–Fejérian analysis of some optimization algorithms. *Studies in Computational Mathematics*, 8:115–152, 2001.
18. P. L. Combettes. Solving monotone inclusions via compositions of nonexpansive averaged operators. *Optimization*, 53(5-6):475–504, 2004.
19. P. L. Combettes and J. C. Pesquet. Primal–dual splitting algorithm for solving inclusions with mixtures of composite, Lipschitzian, and parallel-sum type monotone operators. *Set-Valued and Variational Analysis*, 20(2):307–330, 2012.
20. P. L. Combettes and V. R. Wajs. Signal recovery by proximal Forward–Backward splitting. *Multiscale Modeling & Simulation*, 4(4):1168–1200, 2005.
21. R. Cominetti, J. A. Soto, and J. Vaisman. On the rate of convergence of Krasnoselski–Mann iterations and their connection with sums of Bernoullis. *Israel Journal of Mathematics*, 199(2):757–772, 2014.
22. L. Condat. A primal–dual splitting method for convex optimization involving Lipschitzian, proximable and linear composite terms. *Journal of Optimization Theory and Applications*, 158(2):460–479, 2013.
23. L. Demanet and X. Zhang. Eventual linear convergence of the Douglas–Rachford iteration for basis pursuit. *Mathematics of Computation*, 2015. DOI http://dx.doi.org/10.1090/mcom/2965
24. A. L. Dontchev, M. Quincampoix, and N. Zlateva. Aubin criterion for metric regularity. *Journal of Convex Analysis*, 13:281–297, 2006.
25. A. L. Dontchev and R. T. Rockafellar. *Implicit functions and solution mappings: A view from variational analysis*. Springer, 2009.
26. J. Douglas and H. H. Rachford. On the numerical solution of heat conduction problems in two and three space variables. *Transactions of the American mathematical Society*, 82(2):421–439, 1956.
27. J. Eckstein and D. P. Bertsekas. On the Douglas–Rachford splitting method and the proximal point algorithm for maximal monotone operators. *Mathematical Programming*, 55(1-3):293–318, 1992.
28. M. Fortin and R. Glowinski. *Augmented Lagrangian methods: applications to the numerical solution of boundary-value problems*. Elsevier, 2000.
29. D. Gabay. Chapter IX: Applications of the method of multipliers to variational inequalities. *Studies in mathematics and its applications*, 15:299–331, 1983.
30. D. Gabay and B. Mercier. A dual algorithm for the solution of nonlinear variational problems via finite element approximation. *Computers & Mathematics with Applications*, 2(1):17–40, 1976.
31. R. Glowinski and P. Le Tallec. *Augmented Lagrangian and operator-splitting methods in nonlinear mechanics*, volume 9. SIAM, 1989.
32. B. He and X. Yuan. On non-ergodic convergence rate of Douglas–Rachford alternating direction method of multipliers. *Numerische Mathematik*, 130(3):567–577, 2012.
33. B. He and X. Yuan. On convergence rate of the Douglas–Rachford operator splitting method. *Mathematical Programming*, 2014. DOI http://dx.doi.org/10.1007/s10107-014-0805-x
34. R. Hesse and D. R. Luke. Nonconvex notions of regularity and convergence of fundamental algorithms for feasibility problems. *SIAM Journal on Optimization*, 23(4):2397–2419, 2013.
35. M. A. Krasnosel'skii. Two remarks on the method of successive approximations. *Uspekhi Matematicheskikh Nauk*, 10(1):123–127, 1955.
36. B. Lemaire. Stability of the iteration method for non expansive mappings. *Serdica Mathematical Journal*, 22(3):331p–340p, 1996.
37. B. Lemaire. Which fixed point does the iteration method select? In *Recent Advances in Optimization*, pages 154–167. Springer, 1997.
38. D. Leventhal. Metric subregularity and the proximal point method. *Journal of Mathematical Analysis and Applications*, 360(2):681–688, 2009.
39. A. Lewis, D. Luke, and J. Malick. Local linear convergence for alternating and averaged nonconvex projections. *Found. Comput. Math.*, 9(4):485–513, 2009.
40. G. Li and B. S. Mordukhovich. Hölder metric subregularity with applications to proximal point method. *SIAM Journal on Optimization*, 22(4):1655–1684, 2012.
41. J. Liang, M.J. Fadili, and G. Peyré. Local linear convergence of Forward–Backward under partial smoothness. In *Advances in Neural Information Processing Systems (NIPS)*, pages 1970–1978, 2014.
42. J. Liang, M.J. Fadili, G. Peyré, and Russell Luke. Activity identification and local linear convergence of Douglas-Rachford/ADMM under partial smoothness. In *SSVM 2015 - International Conference on Scale Space and Variational Methods in Computer Vision*, 2015.
43. P. L. Lions and B. Mercier. Splitting algorithms for the sum of two nonlinear operators. *SIAM Journal on Numerical Analysis*, 16(6):964–979, 1979.
44. W. R. Mann. Mean value methods in iteration. *Proceedings of the American Mathematical Society*, 4(3):506–510, 1953.

45. B. Martinet. Brève communication. régularisation d'inéquations variationnelles par approximations successives. *ESAIM: Mathematical Modelling and Numerical Analysis-Modélisation Mathématique et Analyse Numérique*, 4(R3):154–158, 1970.
46. B. Mercier. Topics in finite element solution of elliptic problems. *Lectures on Mathematics*, 63, 1979.
47. G. J. Minty. Monotone (nonlinear) operators in Hilbert space. *Duke Mathematical Journal*, 29(3):341–346, 1962.
48. R. DC Monteiro and B. F. Svaiter. On the complexity of the hybrid proximal extragradient method for the iterates and the ergodic mean. *SIAM Journal on Optimization*, 20(6):2755–2787, 2010.
49. J. J. Moreau. Décomposition orthogonale d'un espace Hilbertien selon deux cônes mutuellement polaires. *CR Acad. Sci. Paris*, 255:238–240, 1962.
50. J. J. Moreau. Proximité et dualité dans un espace Hilbertien. *Bulletin de la Société mathématique de France*, 93:273–299, 1965.
51. Y. Nesterov. *Introductory lectures on convex optimization: A basic course*, volume 87. Springer, 2004.
52. N. Ogura and I. Yamada. Non-strictly convex minimization over the fixed point set of an asymptotically shrinking nonexpansive mapping. *Numerical Functional Analysis and Optimization*, 23(1-2):113–137, 2002.
53. Z. Opial. Weak convergence of the sequence of successive approximations for nonexpansive mappings. *Bull. Amer. Math. Soc*, 73(4):591–597, 1967.
54. G. B. Passty. Ergodic convergence to a zero of the sum of monotone operators in Hilbert space. *Journal of Mathematical Analysis and Applications*, 72(2):383–390, 1979.
55. D. W. Peaceman and H. H. Rachford, Jr. The numerical solution of parabolic and elliptic differential equations. *Journal of the Society for Industrial & Applied Mathematics*, 3(1):28–41, 1955.
56. H. M. Phan. Linear convergence of the Douglas–Rachford method for two closed sets. Technical Report arXiv:1401.6509v1, 2014.
57. B. T. Polyak. *Introduction to optimization*. Optimization Software, 1987.
58. H. Raguet, J. Fadili, and G. Peyré. A generalized Forward–Backward splitting. *SIAM Journal on Imaging Sciences*, 6(3):1199–1226, 2013.
59. S. Reich. Weak convergence theorems for nonexpansive mappings in Banach spaces. *J. Math. Anal. Appl.*, 67:274–276, 1979.
60. R. T. Rockafellar. Monotone operators and the proximal point algorithm. *SIAM Journal on Control and Optimization*, 14(5):877–898, 1976.
61. M. V. Solodov. A class of decomposition methods for convex optimization and monotone variational inclusions via the hybrid inexact proximal point framework. *Optimization Methods and Software*, 19(5):557–575, 2004.
62. M. V. Solodov and B. F. Svaiter. A hybrid approximate extragradient–proximal point algorithm using the enlargement of a maximal monotone operator. *Set-Valued Analysis*, 7(4):323–345, 1999.
63. B. F. Svaiter. A class of Fejer convergent algorithms, approximate resolvents and the hybrid proximal-extragradient method. *Journal of Optimization Theory and Applications*, 162(1):133–153, 2014.
64. P. Tseng. Alternating projection-proximal methods for convex programming and variational inequalities. *SIAM Journal on Optimization*, 7(4):951–965, 1997.
65. B. C. Vũ. A splitting algorithm for dual monotone inclusions involving cocoercive operators. *Advances in Computational Mathematics*, 38(3):667–681, 2013.

# Supplementary Material

**Paper title** Convergence Rates with Inexact Non-expansive Operators
**Journal name** Mathematical Programming series A
**Authors** Jingwei Liang*, Jalal Fadili* and Gabriel Peyré†
*GREYC, CNRS-ENSICAEN-Université de Caen, Jalal.Fadili@ensicaen.fr
†CNRS, CEREMADE, Université Paris-Dauphine

## A Introduction

In this supplementary material, we first present an extra iteration-complexity bound for GFB, and discuss in detail on Vũ's Primal–Dual splitting method. Then to illustrate the usefulness of our results, we apply them to a large class of convex optimization problems and the induced structured monotone inclusions. Experiments on some large scale inverse problems in signal and image processing are shown.

## B GFB method

Since $\boldsymbol{T}_{1,\gamma}$ is firmly non-expansive, then owing to [7, Proposition 4.15], there exists a maximal monotone operator $\boldsymbol{A}_\gamma : \boldsymbol{\mathcal{H}} \to \boldsymbol{\mathcal{H}}$ such that $\boldsymbol{T}_{1,\gamma} = J_{\boldsymbol{A}_\gamma} = (\mathbf{Id} + \boldsymbol{A}_\gamma)^{-1}$. Therefore, for the fixed-point iteration, the following equivalence holds

$$\boldsymbol{z}^\star \in \mathrm{fix}\boldsymbol{T}_\gamma \Leftrightarrow \boldsymbol{z}^\star \in \mathrm{zer}\big(\tfrac{1}{\gamma}\boldsymbol{A}_\gamma + \boldsymbol{B}_{\boldsymbol{\mathcal{S}}}\big). \tag{B.1}$$

We can thus establish a iteration-complexity bounds of GFB for (B.1).

Denote $\boldsymbol{\varepsilon}_k = (\varepsilon_{i,k})_i$, and define $\boldsymbol{p}_{k+1} = \frac{1}{\lambda_k}(\boldsymbol{z}_{k+1} - (1-\lambda_k)\boldsymbol{z}_k) - \boldsymbol{\varepsilon}_k$, $\boldsymbol{g}_{k+1} = \frac{1}{\gamma}\boldsymbol{z}_k - \boldsymbol{B}_{\boldsymbol{\mathcal{S}}}\boldsymbol{z}_k - \frac{1}{\gamma}\boldsymbol{p}_{k+1}$, and $\bar{\boldsymbol{p}}_{k+1} = \frac{1}{k+1}\sum_{j=0}^{k}\boldsymbol{p}_{j+1}$, $\bar{\boldsymbol{g}}_{k+1} = \frac{1}{\gamma}\bar{\boldsymbol{z}}_k - \boldsymbol{B}_{\boldsymbol{\mathcal{S}}}\bar{\boldsymbol{z}}_k - \frac{1}{\gamma}\bar{\boldsymbol{p}}_{k+1}$. Also for the sake of convenience, we assume sequence $(\lambda_k)_{k\in\mathbb{N}} \in [\frac{1}{2\alpha}, \frac{1}{\alpha}[$ is non-decreasing.

**Proposition B.1.** *We have $\boldsymbol{g}_{k+1} \in \frac{1}{\gamma}\boldsymbol{A}_\gamma\boldsymbol{p}_{k+1}$, and under the conditions of Corollary 1 and Theorem 2, there holds,*

$$d\big(0, \tfrac{1}{\gamma}\boldsymbol{A}_\gamma\boldsymbol{p}_{k+1} + \boldsymbol{B}_{\boldsymbol{\mathcal{S}}}\boldsymbol{p}_{k+1}\big) \le \frac{1}{\gamma}\sqrt{\frac{d_0^2 + C_1}{\underline{\tau}_k(k+1)}}, \quad \|\bar{\boldsymbol{g}}_{k+1} + \boldsymbol{B}_{\boldsymbol{\mathcal{S}}}\bar{\boldsymbol{p}}_{k+1}\| \le \frac{2(d_0 + C_2)}{\gamma\lambda_0(k+1)},$$

*where $C_1$ and $C_2$ is similar to those in Corollary 1 and Theorem 2.*

**Proof.** From the fixed-point iteration (17), we have

$$\begin{aligned} &\boldsymbol{z}_{k+1} = \boldsymbol{z}_k + \lambda_k\big((\mathbf{Id} + \boldsymbol{A}_\gamma)^{-1}(\mathbf{Id} - \gamma\boldsymbol{B}_{\boldsymbol{\mathcal{S}}})\boldsymbol{z}_k + \boldsymbol{\varepsilon}_k - \boldsymbol{z}_k\big) \\ \Leftrightarrow\ &\frac{\boldsymbol{z}_{k+1} - (1-\lambda_k)\boldsymbol{z}_k}{\lambda_k} - \boldsymbol{\varepsilon}_k = (\mathbf{Id} + \boldsymbol{A}_\gamma)^{-1}(\mathbf{Id} - \gamma\boldsymbol{B}_{\boldsymbol{\mathcal{S}}})\boldsymbol{z}_k \\ \Leftrightarrow\ &\big(\tfrac{1}{\gamma}\boldsymbol{z}_k - \boldsymbol{B}_{\boldsymbol{\mathcal{S}}}\boldsymbol{z}_k\big) - \big(\tfrac{1}{\gamma}\boldsymbol{p}_{k+1} - \boldsymbol{B}_{\boldsymbol{\mathcal{S}}}\boldsymbol{p}_{k+1}\big) \in \big(\tfrac{1}{\gamma}\boldsymbol{A}_\gamma + \boldsymbol{B}_{\boldsymbol{\mathcal{S}}}\big)\boldsymbol{p}_{k+1}. \end{aligned}$$

By Baillon–Haddad theorem [1], $\mathbf{Id} - \gamma\boldsymbol{B}_{\boldsymbol{\mathcal{S}}} \in \mathcal{A}(\frac{\gamma}{2\beta})$ is non-expansive, then

$$\begin{aligned} d\big(0, \tfrac{1}{\gamma}\boldsymbol{A}_\gamma\boldsymbol{p}_{k+1} + \boldsymbol{B}_{\boldsymbol{\mathcal{S}}}\boldsymbol{p}_{k+1}\big) &\le \|\boldsymbol{g}_{k+1} + \boldsymbol{B}_{\boldsymbol{\mathcal{S}}}\boldsymbol{p}_{k+1}\| = \frac{1}{\gamma}\|(\mathbf{Id} - \gamma\boldsymbol{B}_{\boldsymbol{\mathcal{S}}})\boldsymbol{z}_k - (\mathbf{Id} - \gamma\boldsymbol{B}_{\boldsymbol{\mathcal{S}}})\boldsymbol{p}_{k+1}\| \\ &\le \frac{1}{\gamma}\|\boldsymbol{z}_k - \boldsymbol{p}_{k+1}\| = \frac{1}{\gamma}\|(\boldsymbol{z}_k - \boldsymbol{z}_{k+1})/\lambda_k + \boldsymbol{\varepsilon}_k\| = \frac{1}{\gamma}\|\boldsymbol{e}_k\| \le \frac{1}{\gamma}\sqrt{\frac{d_0^2 + C_1}{\underline{\tau}_k(k+1)}}. \end{aligned}$$

Since we assume $(\lambda_k)_{k\in\mathbb{N}} \in [\frac{1}{2\alpha}, \frac{1}{\alpha}[$ is non-decreasing, hence $\frac{1}{\lambda_0} \geq \frac{1}{\lambda_j}$, $j = 0, ..., k$, then from Theorem 2,

$$\begin{aligned}\|\bar{\boldsymbol{g}}_{k+1} + \boldsymbol{B}_{\boldsymbol{S}}\bar{\boldsymbol{p}}_{k+1}\| &= \frac{1}{\gamma}\|(\mathbf{Id} - \gamma\boldsymbol{B}_{\boldsymbol{S}})\bar{\boldsymbol{z}}_k - (\mathbf{Id} - \gamma\boldsymbol{B}_{\boldsymbol{S}})\bar{\boldsymbol{p}}_{k+1}\| \\ &\leq \frac{1}{\gamma}\|\bar{\boldsymbol{z}}_k - \bar{\boldsymbol{p}}_{k+1}\| = \frac{1}{\gamma(k+1)}\|\textstyle\sum_{j=0}^{k}\frac{1}{\lambda_j}(\boldsymbol{z}_j - \boldsymbol{z}_{j+1} + \lambda_j\boldsymbol{\varepsilon}_j)\| \\ &\leq \frac{1}{\gamma\lambda_0(k+1)}\big(\|\boldsymbol{z}_0 - \boldsymbol{z}_{k+1}\| + \textstyle\sum_{j=0}^{k}\lambda_j\|\boldsymbol{\varepsilon}_j\|\big) \leq \frac{2(d_0 + C_2)}{\gamma\lambda_0(k+1)}.\end{aligned}$$

□

# C Vũ's Primal–Dual splitting method

In this section, we present the iteration-complexity bounds the Primal–Dual splitting method [10]. Define $\boldsymbol{t}_{k+1} = \boldsymbol{T}\boldsymbol{z}_k$ and

$$\boldsymbol{e}_k = (\mathbf{Id} - \boldsymbol{T})\boldsymbol{z}_k = \boldsymbol{z}_k - \boldsymbol{t}_{k+1}.$$

Let $\underline{\tau} = \inf_{k\in\mathbb{N}} \lambda_k(\frac{4\eta\beta-1}{2\eta\beta} - \lambda_k)$, $\overline{\tau} = \sup_{k\in\mathbb{N}} \lambda_k(\frac{4\eta\beta-1}{2\eta\beta} - \lambda_k)$, $\delta = \max\{\frac{1}{\tau}, \frac{1}{\sigma_1}, ..., \frac{1}{\sigma_n}\}$. Define the distance $d_0 = \inf_{\boldsymbol{z}^\star\in\mathrm{fix}\boldsymbol{T}} \|\boldsymbol{z}_0 - \boldsymbol{z}^\star\|_{\boldsymbol{\mathcal{K}}}, \nu_1 = 2\sup_{k\in\mathbb{N}} \|\boldsymbol{T}\boldsymbol{z}_k - \boldsymbol{z}^\star\|_{\boldsymbol{\mathcal{K}}} + \sup_{k\in\mathbb{N}} \lambda_k\|\boldsymbol{\varepsilon}_k\|_{\boldsymbol{\mathcal{K}}}$ and $\nu_2 = 2\sup_{k\in\mathbb{N}} \|\boldsymbol{e}_k - \boldsymbol{e}_{k+1}\|_{\boldsymbol{\mathcal{K}}}$.

Define $\boldsymbol{v}_k = \boldsymbol{z}_k - \boldsymbol{z}_{k+1} + \lambda_k\boldsymbol{\varepsilon}_k$, $\boldsymbol{w}_{k+1} = (\boldsymbol{z}_{k+1} - (1-\lambda_k)\boldsymbol{z}_k)/\lambda_k - \boldsymbol{\varepsilon}_k$, $\boldsymbol{g}_{\boldsymbol{F},k+1} = \boldsymbol{z}_k - \boldsymbol{B}\boldsymbol{z}_k - \boldsymbol{w}_{k+1}$.

**Proposition C.1.** *For the Primal–Dual splitting method [10], there holds*

$$\|\boldsymbol{e}_k\|_{\boldsymbol{\mathcal{K}}} \leq \frac{2\delta}{\eta}\sqrt{\frac{d_0^2 + C_1}{\underline{\tau}(k+1)}} \quad \text{and} \quad \|\boldsymbol{v}_k\|_{\boldsymbol{\mathcal{K}}} \leq \frac{2\delta\lambda_k}{\eta}\sqrt{\frac{d_0^2 + C_1}{\underline{\tau}(k+1)}},$$

*where* $C_1 = \nu_1 \sum_{j\in\mathbb{N}} \lambda_j\|\boldsymbol{\varepsilon}_j\|_{\boldsymbol{\mathcal{K}}} + \nu_2\overline{\tau}\sum_{\ell\in\mathbb{N}} (\ell+1)\|\boldsymbol{\varepsilon}_\ell\|_{\boldsymbol{\mathcal{K}}} < +\infty$.

**Proof.** First, we show that $\|\cdot\|_{\boldsymbol{F}}$ is lower-/upper-bounded by $\|\cdot\|_{\boldsymbol{\mathcal{K}}}$, define operator $\boldsymbol{V}$

$$\boldsymbol{V} : \mathcal{H} \to \boldsymbol{\mathcal{G}},\ x \mapsto \big(\sqrt{\sigma_1}^{-1}L_1x, \cdots, \sqrt{\sigma_n}^{-1}L_nx\big)$$

then, for $\forall x \in \mathcal{H}$,

$$\|\boldsymbol{V}x\|_{\boldsymbol{\mathcal{G}}}^2 = \textstyle\sum_{i=1}^{n}\omega_i\sigma_i\|L_ix\|^2 \leq \|x\|^2\sum_{i=1}^{n}\omega_i\sigma_i\|L_i\|^2 \implies \|\boldsymbol{V}\|^2 \leq \sum_{i=1}^{n}\omega_i\sigma_i\|L_i\|^2,$$

set $\theta = \big(\tau\sum_{i=1}^{n}\omega_i\sigma_i\|L_i\|^2\big)^{-1/2} - 1$, then $\theta > 0$ and

$$\tau(1+\theta)\|\boldsymbol{V}\|^2 \leq \tau(1+\theta)\textstyle\sum_{i=1}^{n}\omega_i\sigma_i\|L_i\|^2 = \frac{1}{1+\theta},$$

since $\tau\sum_i\omega_i\sigma_i\|L_i\|^2 < 1$. Then for $\forall\boldsymbol{z} \in \boldsymbol{\mathcal{K}}$, we have

$$\begin{aligned}\langle\boldsymbol{z}, \boldsymbol{F}\boldsymbol{z}\rangle_{\boldsymbol{\mathcal{K}}} &= \tau^{-1}\|x\|^2 + \textstyle\sum_{i=1}^{n}\sigma_i^{-1}\omega_i\|v_i\|^2 - 2\sum_{i=1}^{n}\omega_i\langle L_ix, v_i\rangle \\ &\leq \tau^{-1}\|x\|^2 + \textstyle\sum_{i=1}^{n}\sigma_i^{-1}\omega_i\|v_i\|^2 + \sum_{i=1}^{n}\omega_i\big(\sigma_i\|L_ix\|^2 + \sigma_i^{-1}\|v_i\|^2\big) \\ &\leq 2\big(\tau^{-1}\|x\|^2 + \textstyle\sum_{i=1}^{n}\sigma_i^{-1}\omega_i\|v_i\|^2\big) \leq 2\delta\big(\|x\|^2 + \sum_{i=1}^{n}\omega_i\|v_i\|^2\big) = 2\delta\|\boldsymbol{z}\|_{\boldsymbol{\mathcal{K}}},\end{aligned}$$

where $\delta = \max\{\tau^{-1}, {\sigma_1}^{-1}, ..., {\sigma_n}^{-1}\}$, and from [10, Equation 3.20], we have

$$\eta\|\boldsymbol{z}\|_{\boldsymbol{\mathcal{K}}}^2 \leq \langle\boldsymbol{z}, \boldsymbol{F}\boldsymbol{z}\rangle_{\boldsymbol{\mathcal{K}}},$$

therefore,

$$\eta\|\boldsymbol{z}\|_{\boldsymbol{\mathcal{K}}}^2 \leq \|\boldsymbol{z}\|_{\boldsymbol{F}}^2 \leq 2\delta\|\boldsymbol{z}\|_{\boldsymbol{\mathcal{K}}}^2. \tag{C.1}$$

Define $\tilde{d}_0 = \inf_{\boldsymbol{z}^\star\in\mathrm{fix}\boldsymbol{T}} \|\boldsymbol{z}_0 - \boldsymbol{z}^\star\|_{\boldsymbol{F}}$, $\tilde{\nu}_1 = 2\sup_{k\in\mathbb{N}} \|\boldsymbol{T}\boldsymbol{z}_k - \boldsymbol{z}^\star\|_{\boldsymbol{F}} + \sup_{k\in\mathbb{N}} \lambda_k\|\boldsymbol{\varepsilon}_k\|_{\boldsymbol{F}}$, $\tilde{\nu}_2 = 2\sup_{k\in\mathbb{N}} \|\boldsymbol{e}_k - \boldsymbol{e}_{k+1}\|_{\boldsymbol{F}}$, $\tilde{C}_1 = \nu_1\sum_{j\in\mathbb{N}} \lambda_j\|\boldsymbol{\varepsilon}_j\|_{\boldsymbol{F}} + \nu_2\overline{\tau}\sum_{\ell\in\mathbb{N}} (\ell+1)\|\boldsymbol{\varepsilon}_\ell\|_{\boldsymbol{F}} < +\infty$. From the bounds (C.1), we have

$$\tilde{d}_0 \leq 2\delta d_0,\ \tilde{\nu}_1 \leq 2\delta\nu_1,\ \text{and}\ \tilde{\nu}_2 \leq 2\delta\nu_2. \tag{C.2}$$

For the fixed-point iteration, under metric $\boldsymbol{F}$, we have

$$\|\boldsymbol{e}_k\|_{\boldsymbol{F}} \leq \sqrt{\frac{\tilde{d}_0^2 + \tilde{C}_1}{\underline{\tau}(k+1)}}, \quad \|\boldsymbol{v}_k\|_{\boldsymbol{F}} \leq \lambda_k\sqrt{\frac{\tilde{d}_0^2 + \tilde{C}_1}{\underline{\tau}(k+1)}}.$$

then combine (C.1) and (C.2) we obtain the desired result. □

Let $\bar{\boldsymbol{v}}_k = \frac{1}{k+1}\sum_{j=0}^{k} \boldsymbol{v}_j$, $\bar{\boldsymbol{z}}_k = \frac{1}{k+1}\sum_{j=0}^{k} \boldsymbol{z}_j$, $\bar{\boldsymbol{w}}_{k+1} = \frac{1}{k+1}\sum_{j=0}^{k} \boldsymbol{w}_j$.

**Proposition C.2.** *The following statements hold, if $C_2 = \sum_{k\in\mathbb{N}} \lambda_k \|\boldsymbol{\varepsilon}_k\|_{\boldsymbol{\mathcal{K}}} < +\infty$,*

$$\|\bar{\boldsymbol{e}}_k\|_{\boldsymbol{\mathcal{K}}} \le \frac{4\delta(d_0 + C_2)}{\eta\Lambda_k}, \ \|\bar{\boldsymbol{v}}_k\|_{\boldsymbol{\mathcal{K}}} \le \frac{4\delta(d_0 + C_2)}{\eta(k+1)}.$$

**Proof.** The result of a combination of (C.1), Theorem 2, Corollary 2 and Proposition B.1. For $\bar{\boldsymbol{e}}_k$, we have

$$\begin{aligned}\|\bar{\boldsymbol{e}}_k\|_{\boldsymbol{\mathcal{K}}} &\le \frac{1}{\eta}\|\bar{\boldsymbol{e}}_k\|_{\boldsymbol{F}} = \frac{1}{\eta}\|\tfrac{1}{\Lambda_k}\textstyle\sum_{j=0}^{k} \lambda_j \boldsymbol{e}_j\|_{\boldsymbol{F}} = \dfrac{1}{\eta\Lambda_k}\|\textstyle\sum_{j=0}^{k} (\boldsymbol{e}_j - \boldsymbol{e}_{j+1}) + \sum_{j=0}^{k}\lambda_j \boldsymbol{\varepsilon}_j\|_{\boldsymbol{F}} \\ &\le \frac{2\delta}{\eta\Lambda_k}\big(\|\boldsymbol{z}_0 - \boldsymbol{z}^\star\|_{\boldsymbol{\mathcal{K}}} + \|\boldsymbol{z}^\star - \boldsymbol{z}_{k+1}\|_{\boldsymbol{\mathcal{K}}} + \textstyle\sum_{j=0}^{k}\lambda_j\|\boldsymbol{\varepsilon}_j\|_{\boldsymbol{\mathcal{K}}}\big) \le \dfrac{4\delta(d_0 + C_2)}{\eta\Lambda_k}.\end{aligned}$$

Replacing $\Lambda_k$ with $k+1$ obtains the result for $\|\bar{\boldsymbol{v}}_k\|_{\boldsymbol{\mathcal{K}}}$. □

Now, if we reformulate the dual inclusion to the following format,

$$\text{Find } \boldsymbol{v} \in \boldsymbol{\mathcal{G}} \text{ s.t. } (\exists x \in \mathcal{H}), \begin{cases} 0 \in Cx + Bx + \sum_{i=1}^{n}\omega_i L_i^* v_i, \\ 0 \in (A_i^{-1} + D_i^{-1})v_i - (L_i x - r_i), i \in [\![1, n]\!]. \end{cases} \tag{C.3}$$

Then for Vũ's algorithm, define

$$s_{k+1} = J_{\tau C}\big(x_k - \tau(\textstyle\sum_i \omega_i L_i^* v_{i,k} + Bx_k)\big), \ \ t_{i,k+1} = J_{\sigma_i A_i^{-1}}\big(v_{i,k} + \sigma_i(L_i y_{i,k+1} - D_i^{-1} v_{i,k} - r_i)\big),$$

we have

$$\begin{aligned}\big(\tfrac{1}{\tau}\mathrm{Id} - B\big)\, x_k - \big(\tfrac{1}{\tau}\mathrm{Id} - B\big)\, s_{k+1} &\in (C + B)s_{k+1} + \textstyle\sum_{i=1}^{n}\omega_i L_i^* v_{i,k}, \\ \big(\tfrac{1}{\sigma_i}\mathrm{Id} - D_i^{-1}\big)v_{i,k} - \big(\tfrac{1}{\sigma_i}\mathrm{Id} - D_i^{-1}\big)t_{i,k+1} &\in (A_i^{-1} + D_i^{-1})t_{i,k+1} - (L_i y_{k+1} - r_i).\end{aligned}$$

Define operator,

$$\boldsymbol{G} : \boldsymbol{\mathcal{K}} \to 2^{\boldsymbol{\mathcal{K}}}, \, (x, \boldsymbol{v}) \mapsto \big((\tfrac{\mathrm{Id}}{\tau} - B)x\big) \times \big((\tfrac{\mathrm{Id}}{\sigma_1} - D_1^{-1})v_1\big) \times \cdots \times \big((\tfrac{\mathrm{Id}}{\sigma_n} - D_n^{-1})v_n\big),$$

then we have

$$\boldsymbol{g}_{k+1} = \boldsymbol{G}\boldsymbol{z}_k - \boldsymbol{G}\boldsymbol{t}_{k+1} \in \begin{pmatrix} (C + B)s_{k+1} + \sum_i \omega_i L_i^* v_{i,k} \\ (A_1^{-1} + D_1^{-1})t_{1,k+1} - (L_1 y_{k+1} - r_1) \\ \vdots \\ (A_n^{-1} + D_n^{-1})t_{n,k+1} - (L_n y_{k+1} - r_n) \end{pmatrix}.$$

Denote the right hand side term of the inclusion as $\boldsymbol{M}(\boldsymbol{v}_k, \boldsymbol{t}_{k+1}, y_{k+1})$

**Proposition C.3.** *For the Vũ's algorithm, there holds for the dual inclusion* (C.3)

$$d\big(0, \boldsymbol{M}(\boldsymbol{v}_k, \boldsymbol{t}_{k+1}, y_{k+1})\big) \le \frac{2\delta^2}{\eta}\sqrt{\frac{d_0^2 + C_1}{\underline{\tau}(k+1)}}.$$

**Proof.** Again, owing to Baillon–Haddad theorem [1], $\mathrm{Id} - \tau B \in \mathcal{A}(\frac{\tau}{2\mu})$ and $\mathrm{Id} - \sigma_i D_i^{-1} \in \mathcal{A}(\frac{\sigma_i}{2\nu_i})$, $i \in [\![1, n]\!]$ are non-expansive, denote them $G_B$ and $G_{D_i}$ respectively, then we have

$$\begin{aligned} &d\big(0, \boldsymbol{M}(\boldsymbol{v}_k, \boldsymbol{t}_{k+1}, y_{k+1})\big) \le \|\boldsymbol{g}_{k+1}\|_{\boldsymbol{\mathcal{K}}} \\ &= \big(\|\tfrac{1}{\tau}\big(G_B x_k - G_B s_{k+1}\big)\|^2 + \textstyle\sum_{i=1}^{n}\omega_i\|\tfrac{1}{\sigma_i}\big(G_{D_i} v_{i,k} - G_{D_i} t_{i,k+1}\big)\|^2\big)^{\frac{1}{2}} \\ &\le \big(\tfrac{1}{\tau^2}\|x_k - s_{k+1}\|^2 + \textstyle\sum_{i=1}^{n}\omega_i \tfrac{1}{\sigma_i^2}\|v_{i,k} - t_{i,k+1}\|^2\big)^{\frac{1}{2}} \le \delta\big(\|x_k - s_{k+1}\|^2 + \textstyle\sum_{i=1}^{n}\omega_i\|v_{i,k} - t_{i,k+1}\|^2\big)^{\frac{1}{2}} \\ &= \delta\|\boldsymbol{z}_k - \boldsymbol{t}_{k+1}\|_{\boldsymbol{\mathcal{K}}} = \delta\|\boldsymbol{e}_k\|_{\boldsymbol{\mathcal{K}}} \le \frac{2\delta^2}{\eta}\sqrt{\frac{d_0^2 + C_1}{\underline{\tau}(k+1)}}. \end{aligned}$$

□

# D Numerical experiments

To demonstrate the established iteration-complexity bounds and local convergence rate. In this section, we take 3 inverse problems as example:

(1) anisotropic total variation (TV) based deconvolution with box constraint;

(2) matrix completion with non-negativity constraints (NMC);

(3) principal component pursuit problem (PCP) with application to video background and foreground decomposition.

All the problems are solved by both GFB and Vũ's Primal–Dual splitting method (PDS).

## D.1 Anisotropic TV deconvolution

Suppose the blurred observation $y \in \mathbb{R}^n$ is the convolution of $x_0 \in \mathbb{R}^n$ and a point spread function (PSF) $h$ contaminated with additive white Gaussian noise $w$, which reads

$$y = \mathcal{M}(x_0) + w,$$

where $\mathcal{M} : \mathbb{R}^n \to \mathbb{R}^n$ is the linear operator associated to $h$. The deconvolution procedure is to provably recover or approximate $x_0$ from $y$, here we consider the anisotropic TV [9] based deconvolution model which is

$$\min_{x \in \mathbb{R}^n} \tfrac{1}{2}\|\mathcal{M}(x) - y\|^2 + \mu\|\nabla x\|_1 + \iota_\Omega(x), \tag{D.1}$$

where $\mu > 0$ is the regularization parameter determined based on the noise level, $\iota_\Omega(\cdot)$ is the indicator function of box constraint, for instance $\Omega = [0, 255]^n$ if $x$ is gray scale image. The problem can be solved by both GFB and PDS methods, where the proximity operator of $\iota_\Omega(\cdot)$ is the projection onto set $\Omega$, and for GFB method, the proximity operator of $\|\nabla \cdot\|_1$ is computed by minimum graph-cut [5, 6].

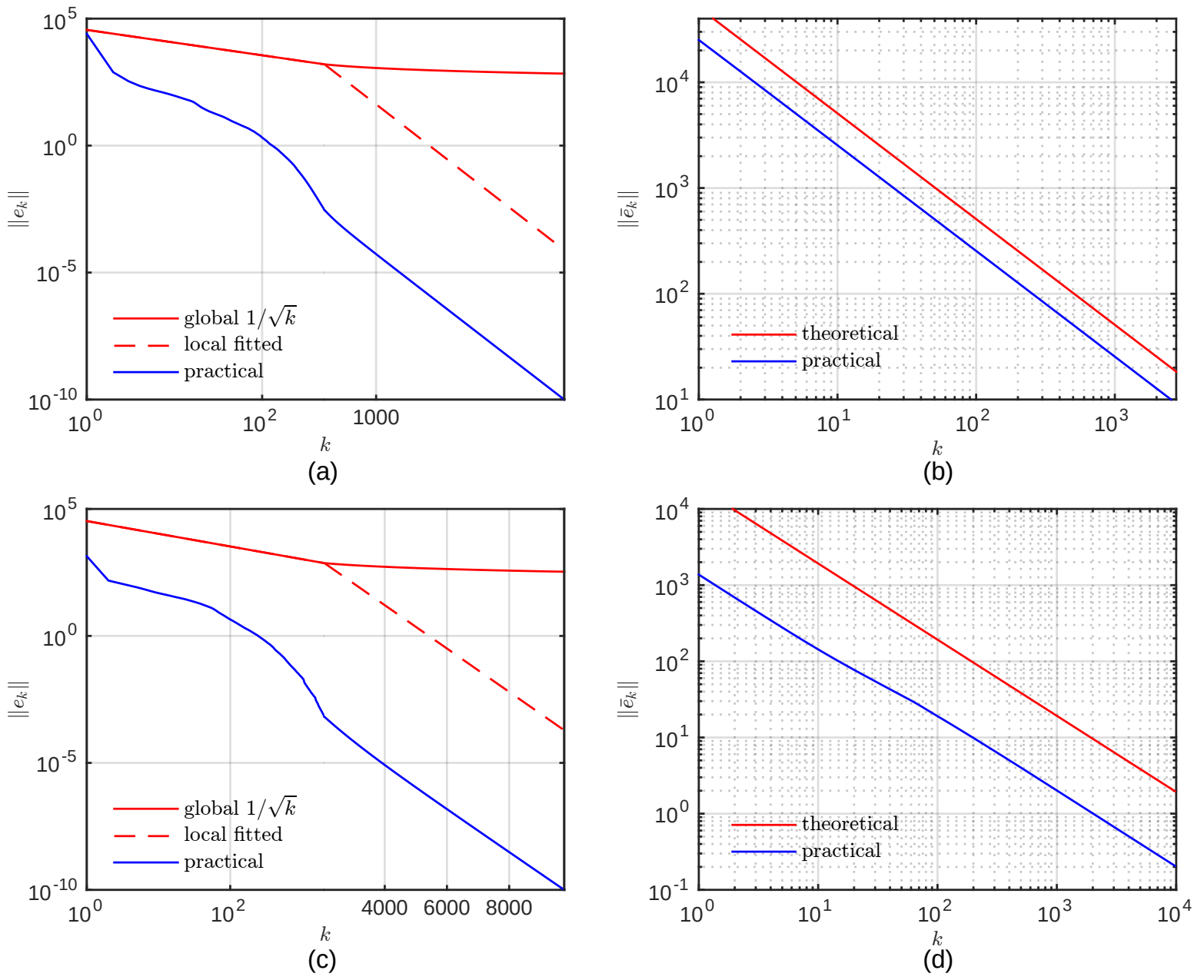

**Figure A:** TV deconvolution of the cameraman image. (a) pointwise convergence of GFB, (b) ergodic convergence of GFB, (c) pointwise convergence of PDS, (d) ergodic convergence of PDS. Both methods achieve local linear convergence. Note that $(\|e_k\|)_{k\in\mathbb{N}}$ is non-increasing, which coincides with Lemma 5 when $\varepsilon_k = 0$.

Figure A displays the observed pointwise and ergodic rates of $\|e_k\|$ and the theoretical bounds given by Theorem 1 and 2. Pointwise convergence rate is shown in subfigure (a) and (c), whose left half is log-log plot while the right half is semi-log plot. As predicted by Theorem 1, globally $\|e_k\|$ converges at the rate of $O\big(1/\sqrt{k}\big)$. Then for a sufficiently large iteration number, a linear convergence regime takes over as clearly seen from the semi-log plot, which is in consistent with the result of Theorem 3. Let us mention that the local linear convergence curve (dashed line) is *fitted* to the observed one, since the regularity modulus necessary to compute the theoretical rate in Theorem 3 is not easy to estimate. For the ergodic convergence, subfigure (b) and (d) of Figure A, $O(1/k)$ convergence rates are observed which coincides with Theorem 2.

## D.2 Non-negative matrix completion

Suppose we observe measurements $y \in \mathbb{R}^p$ of a low rank matrix $X_0 \in \mathbb{R}^{m\times n}$ with non-negative entries

$$y = \mathcal{M}(X_0) + w,$$

where $\mathcal{M} : \mathbb{R}^{m\times n} \to \mathbb{R}^p$ is a measurement operator, and $w$ is the noise. In our experiment here, $\mathcal{M}$ selects $p$ entries of its argument uniformly at random. The matrix completion problem consists in recovering $X_0$, or finding an approximation of it, by solving a convex optimization problem, namely the minimization of the nuclear norm [3, 4, 8]. In penalized form, the problem reads

$$\min_{X\in\mathbb{R}^{m\times n}} \tfrac{1}{2}\|y - \mathcal{M}(X)\|^2 + \mu\|X\|_* + \iota_{P_+}(X), \tag{D.2}$$

where $\iota_{P_+}(\cdot)$ is the indicator function of the non-negative orthant accounting for the non-negativity constraint, and $\mu > 0$ is a regularization parameter typically chosen proportional to the noise level. The proximity operator of both $\|\cdot\|_*$, $\iota_{P_+}(\cdot)$ have explicit forms, since $\mathrm{prox}_{\|\cdot\|_*}(X)$ amounts to soft–thresholding the singular values of $X$ and $\mathrm{prox}_{\iota_{P_+}}(X)$ is the projector on the non-negative orthant.

Figure B displays the observed pointwise and ergodic rates of $\|e_k\|$ and the theoretical bounds computed given by Theorem 1 and 2. Both global and local convergence behaviours are similar to those observed in Figure A.

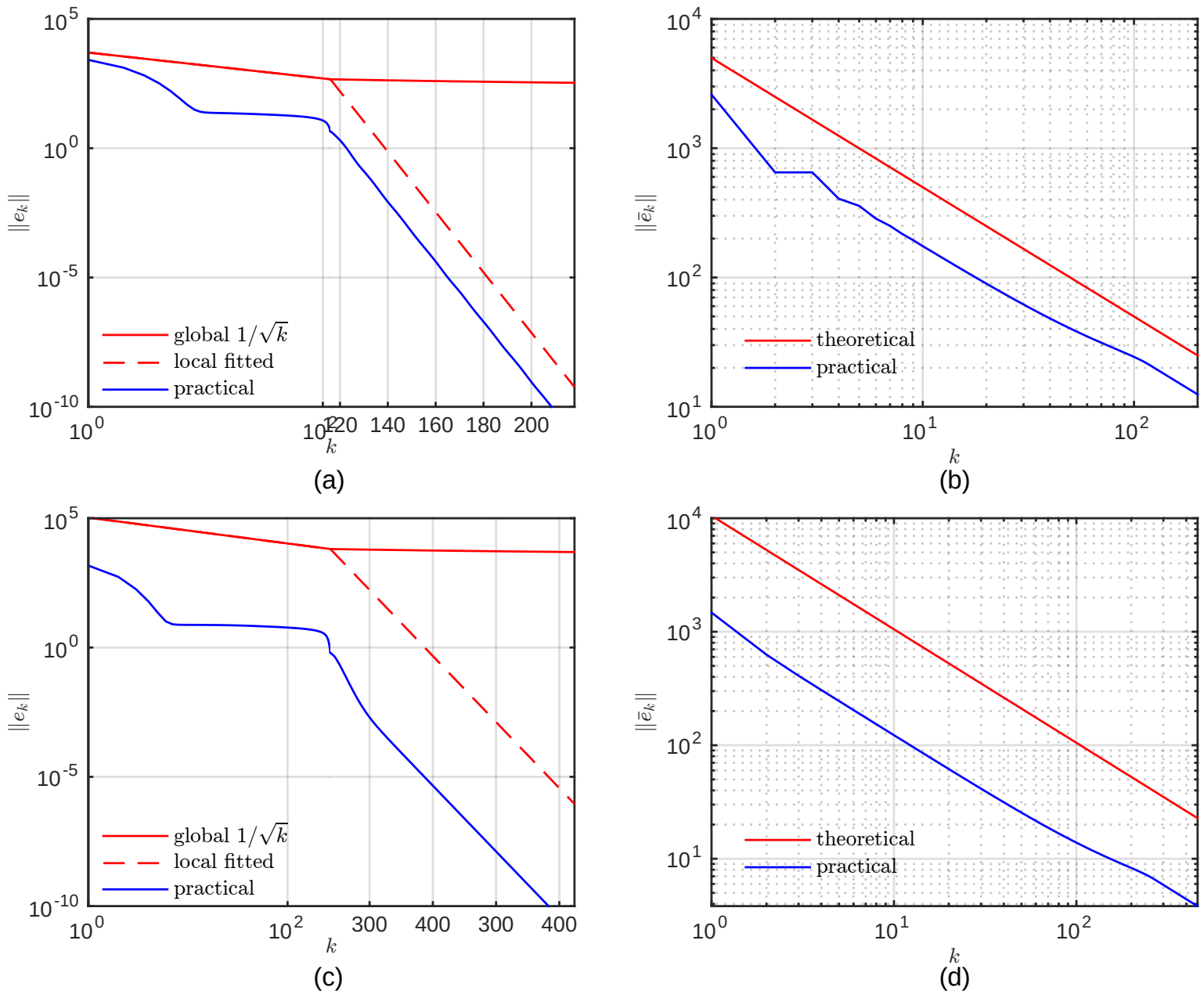

**Figure B:** The size of the matrix $X$ is $400 \times 300$, $\mathrm{rank}(X) = 20$ and the operator $\mathcal{M}$ is random projection mask. (a) pointwise convergence of GFB, (b) ergodic convergence of GFB, (c) pointwise convergence of PDS, (d) ergodic convergence of PDS.

## D.3 Principal component pursuit

In this experiment, we consider the PCP problem [2], and apply it to decompose a video sequence into its background and foreground components. The rationale behind this is that since the background is virtually the same in all frames, if the latter are stacked as columns of a matrix, it is likely to be low–rank (even of rank 1 for perfectly constant background). On the other hand, moving objects appear occasionally on each frame and occupy only a small fraction of it. Thus the corresponding component would be sparse.

Assume that a real matrix $M \in \mathbb{R}^{m\times n}$ can be written as

$$M = X_{L,0} + X_{S,0} + w,$$

where $X_{L,0}$ is low–rank, $X_{S,0}$ is sparse and $w$ is a perturbation matrix of variance $\sigma$ that accounts for model imperfection (noise). The PCP proposed in [2] attempts to provably recover $(X_{L,0}, X_{S,0})$ to a good approximation, by solving a convex optimization. Here, toward an application to video decomposition, we also add a non-negativity constraint to the low–rank component, which leads to the following convex problem

$$\min_{X_L, X_S \in \mathbb{R}^{m\times n}} \tfrac{1}{2}\|M - X_L - X_S\|_F^2 + \mu_1\|X_S\|_1 + \mu_2\|X_L\|_* + \iota_{P_+}(X_L), \tag{D.3}$$

where $\|\cdot\|_F$ is the Frobenius norm.

Observe that for fixed $X_L$, the minimizer of (D.3) is $X_S^\star = \operatorname{prox}_{\mu_1\|\cdot\|_1}(M - X_L)$. Thus, (D.3) is equivalent to

$$\min_{X_L \in \mathbb{R}^{m\times n}} {}^{1}\big(\mu_1\|\cdot\|_1\big)(M - X_L) + \mu_2\|X_L\|_* + \iota_{P_+}(X_L), \tag{D.4}$$

where ${}^{1}\big(\mu_1\|\cdot\|_1\big)(M - X_L) = \min_Z \tfrac{1}{2}\|M - X_L - Z\|_F^2 + \mu_1\|Z\|_1$ is the Moreau Envelope of $\mu_1\|\cdot\|_1$ of index 1, and hence has 1–Lipschitz continuous gradient.

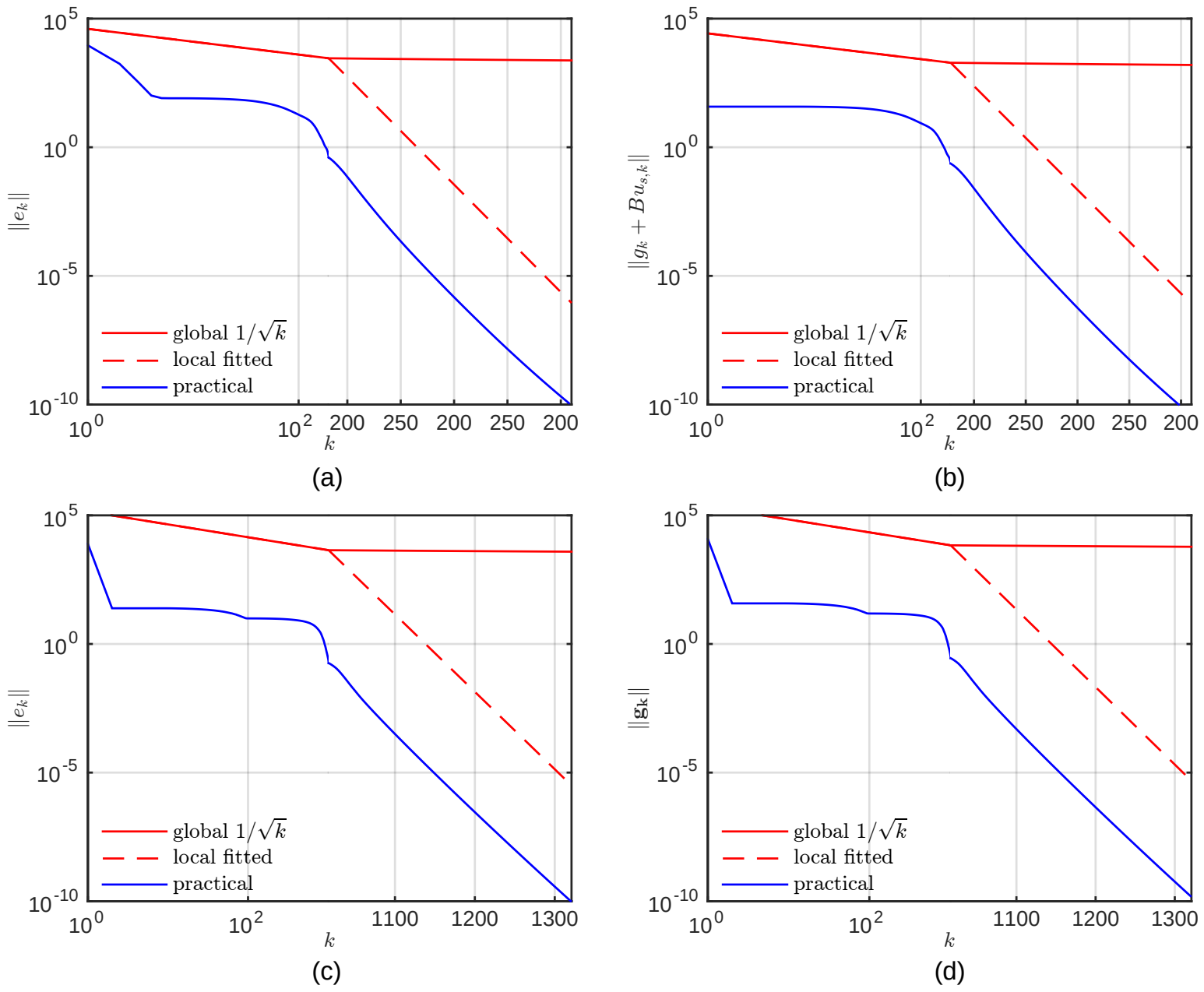

**Figure C:** The size of the matrix $M$ is $400 \times 300$, $\operatorname{rank}(X_{L,0}) = 20$ and the sparsity of $X_{S,0}$ is 25%, namely, 25% of the elements of $X_{L,0}$ are non-zero. (a) pointwise convergence of $\|e_k\|$ of GFB, (b) pointwise convergence of $\|g_k + B(\sum_i \omega_i u_{i,k})\|$ of GFB, (c) pointwise convergence of $\|e_k\|$ of PDS, (d) pointwise convergence $\|\boldsymbol{g}_k\|$ of PDS.

We first use a synthetic example to demonstrate the comparison of the two methods, as shown in Figure C. Pointwise convergence rate of $\|e_k\|$ is shown in subfigure (a) and (c). Then in subfigure (b) and (d), we display the convergence behaviour of the criteria provided in Proposition 4 and C.3.

**Video sequence** The video sequence consists of 400 frames, each of resolution $288 \times 384$ stacked as a column of the matrix $M$. Hence $M$ is of size $110592 \times 400$. We then solved (D.4) to decompose the video into its foreground and background.

For the video test, we choose GFB method to solve the problem since it's faster than PDS. Figure D displays the observed pointwise and ergodic rates and those predicted by Proposition 6 and 7. Figure E shows the decomposition of the video sequence, column (a) shows 3 frames from the video, column (b) and (c) are the corresponding low–rank component $X_L$ and sparse component $X_S$. Notice that $X_L$ correctly recovers the background, while $X_S$ correctly identifies the moving pedestrians and their shadows.

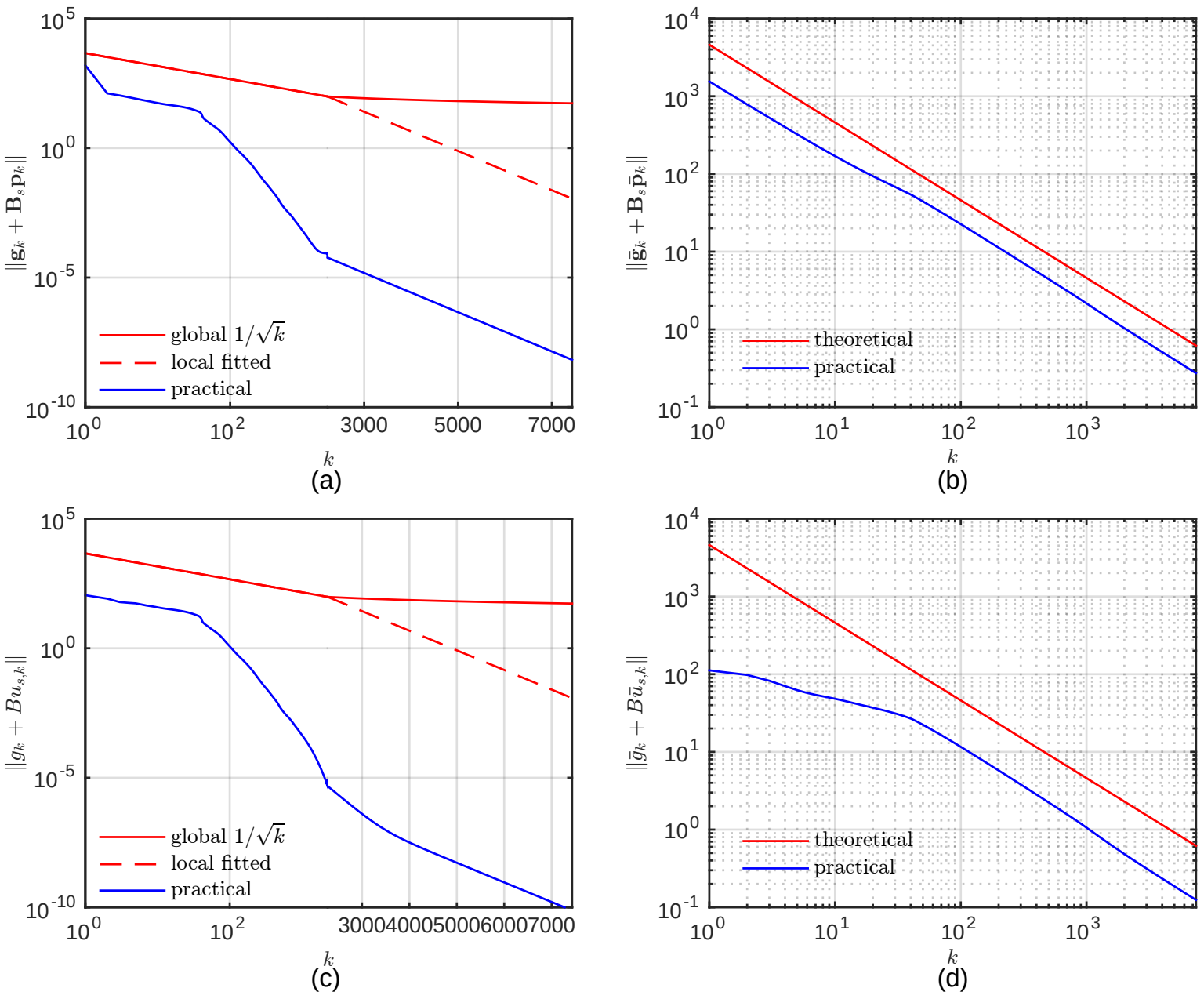

**Figure D:** Observed rates and theoretical bounds of GFB. (a) pointwise convergence of $\|\boldsymbol{g}_k + \boldsymbol{B}_{\boldsymbol{S}}(\boldsymbol{p}_k)\|$, (b) ergodic convergence of $\|\bar{\boldsymbol{g}}_k + \boldsymbol{B}_{\boldsymbol{S}}\bar{\boldsymbol{p}}_k\|$, (c) pointwise convergence of $\|g_k + B(\sum_i \omega_i u_{i,k})\|$, (d) ergodic convergence of $\|\bar{g}_k + B(\sum_i \omega_i \bar{u}_{i,k})\|$.

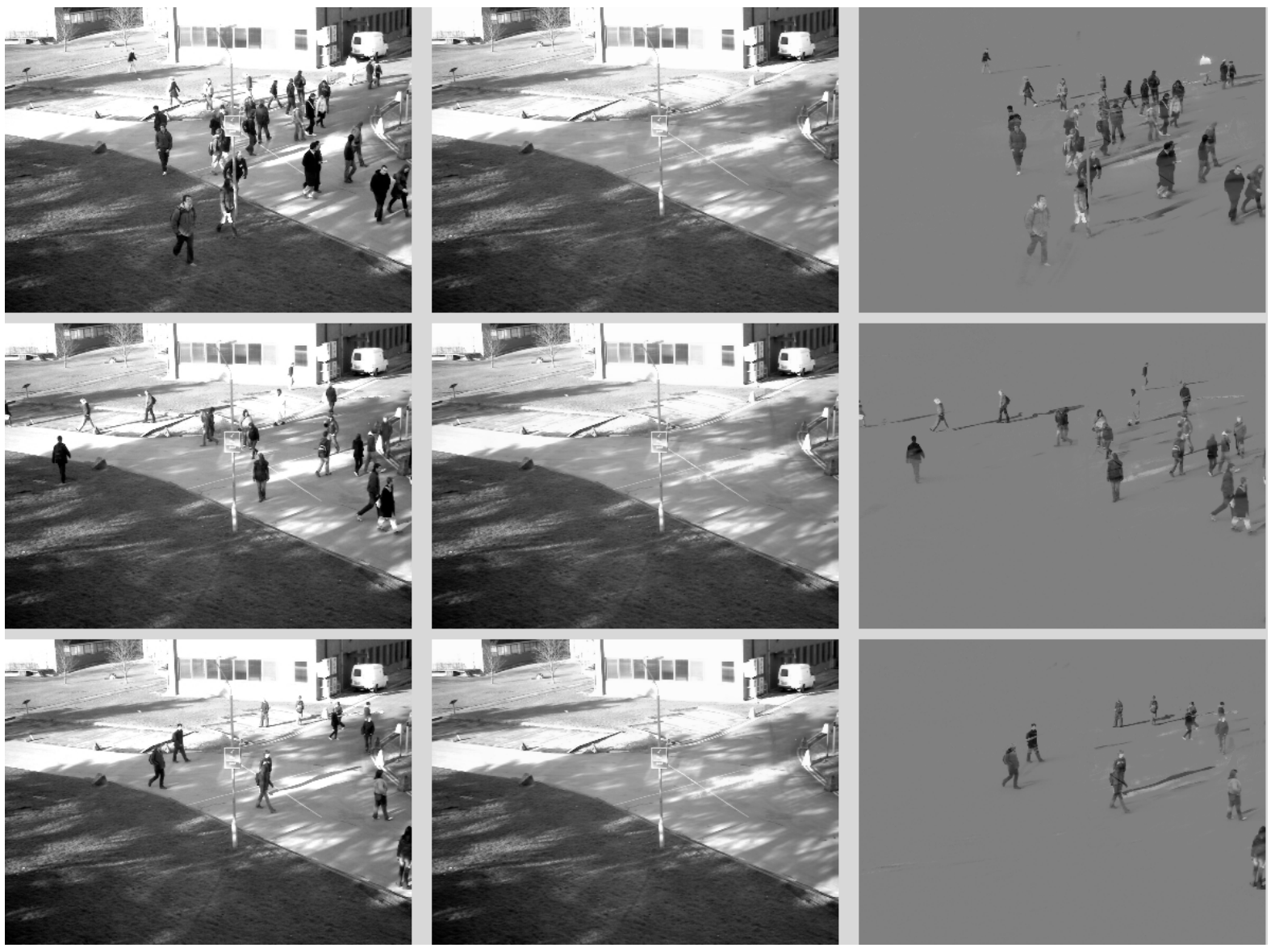

**Figure E:** Three frames from a 400-frame video sequence taken in campus. Left column, original video $M$; Middle column, low rank $X_L$; Right column, sparse $X_S$ obtained by (D.3).

### D.4 Non-stationary iteration

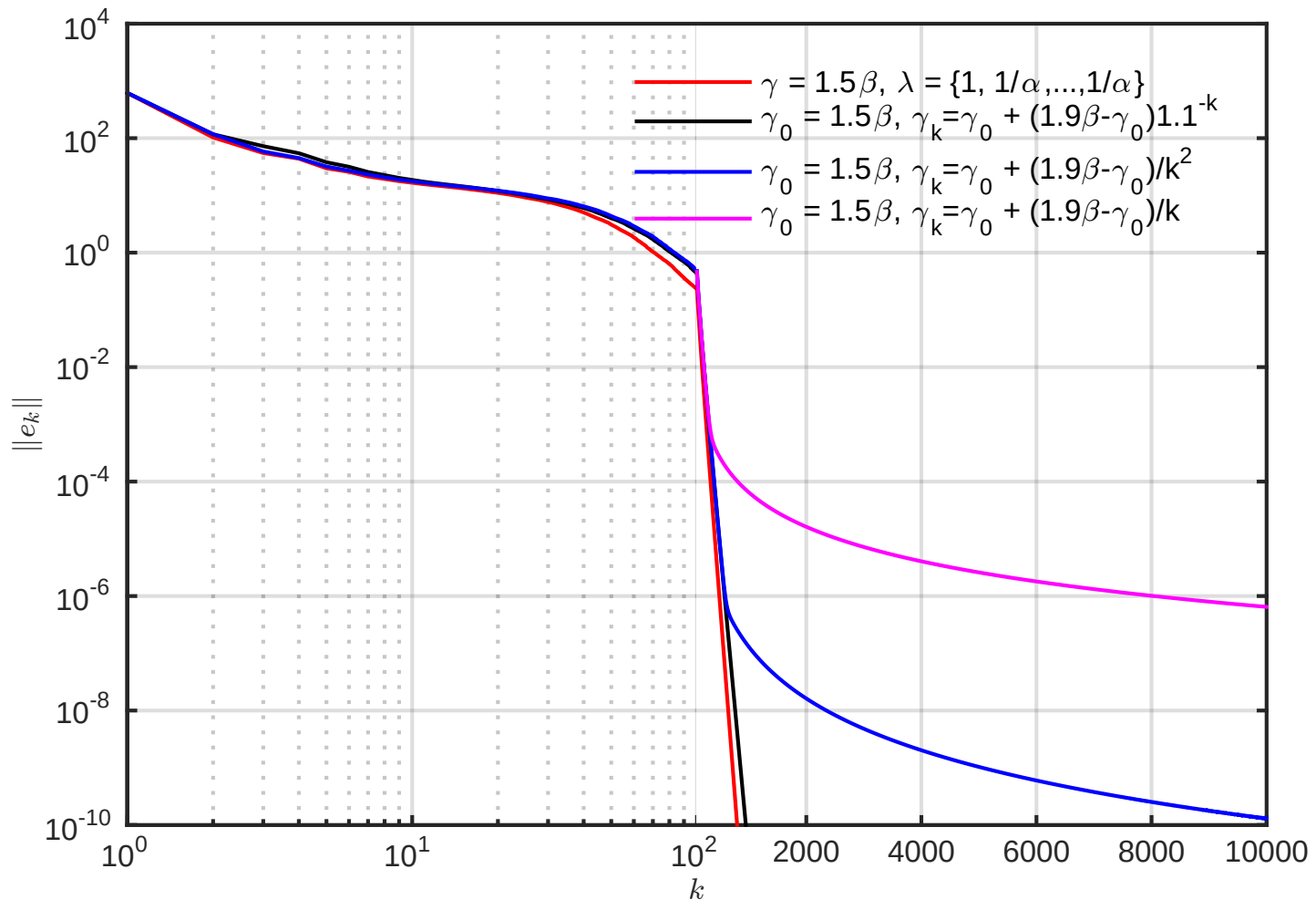

**Figure F:** Comparison of stationary and non-stationary iterations of GFB.

Now we illustrate the non-stationary iteration of GFB applied to PCP problem. The above comparison indicates that in practice we can choose the relaxation parameter $\lambda_k$ as $\lambda_1 = 1$ and $\lambda_k \approx \frac{1}{\alpha}$, $\forall k \geq 2$. Next, we compare this setting with the non-stationary GFB, for the stationary case, we let $\gamma = 1.5\beta, \lambda_k = \{1, \frac{1}{1.05\alpha}, \cdots\}$, for the non-stationary case, let $\gamma_0 = 1.5\beta$, then $3$ scenarios of $\gamma_k$ are considered,

$$\gamma_{1,k} = \gamma_0 + \frac{1.9\beta - \gamma_0}{1.1^k}, \quad \gamma_{2,k} = \gamma_0 + \frac{1.9\beta - \gamma_0}{k^2}, \quad \gamma_{3,k} = \gamma_0 + \frac{1.9\beta - \gamma_0}{k}.$$

The result is given in Figure F. We can conclude from the observation that

- Globally, the non-stationary iterations are slower than the stationary one. Then locally, the local linear rates are also slower than the stationary one, moreover, for $\gamma_{2,k}$ and $\gamma_{3,k}$, the iterations eventually become sublinear, implying that, the iteration is controlled by the perturbation error;
- The non-summability of $(|\gamma_{3,k} - \gamma_0|)_{k\in\mathbb{N}} \notin \ell^1_+$ shows that our assumption on the convergence of the non-stationary iteration is appropriate.

# References

[1] J. B. Baillon and G. Haddad. Quelques propriétés des opérateurs angle-bornés etn-cycliquement monotones. *Israel Journal of Mathematics*, 26(2):137–150, 1977.

[2] E. J. Candès, X. Li, Y. Ma, and J. Wright. Robust principal component analysis? *Journal of the ACM (JACM)*, 58(3):11, 2011.

[3] E. J. Candès and B. Recht. Exact matrix completion via convex optimization. *Foundations of Computational Mathematics*, May 2009.

[4] E. J. Candès and T. Tao. The power of convex relaxation: near-optimal matrix completion. *IEEE Trans. Inform. Theory*, 56(5):2053–2080, 2010.

[5] A. Chambolle and J. Darbon. On total variation minimization and surface evolution using parametric maximum flows. *International journal of computer vision*, 84(3):288–307, 2009.

[6] D. S. Hochbaum. An efficient algorithm for image segmentation, markov random fields and related problems. *Journal of the ACM (JACM)*, 48(4):686–701, 2001.

[7] H. Raguet, J. Fadili, and G. Peyré. Generalized forward–backward splitting. *SIAM Journal on Imaging Sciences*, 6(3):1199–1226, 2013.

[8] B. Recht, M. Fazel, and P. A. Parrilo. Guaranteed minimum-rank solutions of linear matrix equations via nuclear norm minimization. *SIAM Review*, 52(3):471–501, 2010.

[9] L. I Rudin, S. Osher, and E. Fatemi. Nonlinear total variation based noise removal algorithms. *Physica D: Nonlinear Phenomena*, 60(1):259–268, 1992.

[10] B. C. Vũ. A splitting algorithm for dual monotone inclusions involving cocoercive operators. *Advances in Computational Mathematics*, pages 1–15, 2011.